\documentclass[12pt,a4paper]{amsart}
\usepackage[utf8x]{inputenc}
\usepackage[polutonikogreek,greek,english]{babel}

\usepackage{color,graphicx,hyperref}
\usepackage{amsmath, mathrsfs, amsfonts, amssymb}

\newcommand{\defi}[1]{\textsf{#1}} 

\newtheorem{theorem}{Theorem}

\begin{document}

\title[Diophantus Revisited]{Diophantus Revisited: On rational surfaces and K3 surfaces in the ``Arithmetica''}
\author{Ren\'e Pannekoek}
\email{pannekoek@gmail.com}
\address{Imperial College London, South Kensington Campus, 180 Queen's Gate, London SW7 2BZ, United Kingdom}

\begin{abstract}
This article aims to show two things: firstly, that certain problems in Diophantus'  {\it Arithmetica} lead to equations defining del Pezzo surfaces or other rational surfaces, while certain others lead to K3 surfaces; secondly, that Diophantus' own solutions to these problems, at least when viewed through a modern lens, imply the existence of par\-am\-et\-riz\-at\-ions of these surfaces, or of parametrizations of rational curves lying on them. 
\end{abstract}
\maketitle

\tableofcontents

\section{Introduction}

The works of the mathematician Diophantus have often struck readers as idiosyncratic. In the introduction \cite[p. 2]{tannery} to his main surviving work, the {\it Arithmetica}, Diophantus himself comments on the novelty of his science\footnote{All translations in the article are mine, unless otherwise stated.}:
\begin{quotation}
Knowing that you, esteemed Dionysius, are serious about learning to solve problems concerning numbers, I have attempted to organize the method (...)~Perhaps the matter seems difficult, because it is not yet well-known, and the souls of beginners are wont to despair of success (...) 
\end{quotation}
Likewise, the 19${}^\textnormal{th}$-century German historian of mathematics Hermann Hankel was right on the money when he wrote (in \cite[p. 157]{hankel})
\begin{quotation}
Amidst this dreary desolation, there suddenly rises a man with youthful buoyancy: Diophantus. Whence? Whither? Who are his predecessors, who his successors? -- we do not know -- everything a riddle.
\end{quotation}
Hankel's quote beautifully conveys the impression that many others since, the present author among them, have formed of Diophantus' work. Against the dark and depressing backdrop of Hellenistic ``number theory'', Diophantus' {\it Arithmetica} is an absolute miracle of originality. His strange but powerful results prefigure many of the advances that have since been made in the study of ``rational points on varieties''.

In this article, I hope to show the reader some of the ways in which the {\it Arithmetica} has been a precursor of modern number theory. In particular, I will show that Diophantus has solved problems that, for us, would be equivalent to determining rational points on del Pezzo surfaces and K3 surfaces. Furthermore, I hope to convince the reader that, while Diophantus routinely contents himself with finding just a {\it single} solution to his equations, if one just follows the {\it methods} that Diophantus expounds, one frequently ends up with actual parametrizations of rational curves on the surfaces under consideration, or indeed even of the surfaces themselves.\footnote{Here and elsewhere, I take a {\it parametrization} of a variety $X$ to mean a dominant rational map $\mathbb{P}^n \dashrightarrow X$, where $n$ is some positive integer, and where the ground field is invariably $\mathbb{Q}$. Note that I thus do not require parametrizations to be birational; in most cases, they will be, but not in all.} 

\subsection{Diophantus: life and works}

Little is known about the life, or even times, of Diophantus. He certainly lived in Alexandria, and probably did so in the 3${}^{\textnormal{rd}}$ century of the current era.\footnote{For the uncertainties surrounding his dates (and ethnicity), see \cite{schappacher}.} His most important surviving work consists of the ten books of the {\it Arithmetica} that have come down to us: six are extant in the original Greek, four survive only in an Arabic translation that came to light as recently as 1968, while the remaining three are presumed to be lost. In the ten surviving books, some $300$ problems in the theory of numbers are first stated and then solved. 

It is not easy to exaggerate the influence of Diophantus' work on mathematicians throughout the past $1500$ years, first in Islamic, and then later in Western mathematics. It was in a copy of Bachet's 1621 edition of the {\it Arithmetica} that Fermat's Last Theorem was first presented to the world, as a note scribbled in the margin. In the same edition the reader encounters, for the first time in history, the duplication procedure for a certain family of elliptic curves, which Bachet merely states as a commentary to one of Diophantus' problems (Problem VI.19). And then there is the famous ``porism'' of Diophantus: this asserts that a difference between two positive rational cubes can be written also as a {\it sum} of positive rational cubes, again a result which nowadays we would explain in terms of elliptic curves. (For more on the history of the reception of the {\it Arithmetica}, see \cite{schappacher_pre} and \cite{schappacher}). It is in honour of Diophantus' achievements that polynomial equations with rational coefficients, of which it is asked to find rational solutions, are still universally known as ``Diophantine equations''.
        
\subsection{The ``Arithmetica''} As mentioned above, the surviving part of the {\it Arithmetica} consists of ten books. We will use the traditional numbering (from I to VI) for the six Greek ones. It is generally thought that Books I--III actually correspond to the first three books, whereas the four ``Arabic'' books are assumed to have formed Books $4$ through $7$, and we will refer to them as such. The Greek Books IV--VI are thus thought to represent a portion of the work that came {\it after} the four Arabic books. 

As Diophantus writes in the concluding section of the introduction, the problems are roughly ordered ``from the simpler to the more involved'' (\textgreek{ἀπὸ ἁπλουστέρων ἐπὶ	σκολιώτερα}, \cite[p. 16]{tannery}), as the following brief selection demonstrates:
\begin{quotation}
\upshape
{\sc Problem I.1.} To divide a given number into two numbers with given difference.
\end{quotation} 
\begin{quotation}
\upshape
{\sc Problem I.2.} To divide a given number into two numbers with given ratio.
\end{quotation} 
\begin{quotation}
\upshape
{\sc Problem II.8.} To divide a given square into two squares.
\end{quotation} 
\begin{quotation}
\upshape
{\sc Problem II.9.} To divide a given number, which is a sum of two squares, into two different squares.
\end{quotation} 
\begin{quotation}
\upshape
{\sc Problem VI.24.} To find a right-angled triangle such that the perimeter is a cube, and the perimeter added to the area gives a square.
\end{quotation} 
In the solutions, Diophantus always chooses specific values for the ``given'' quantities. So for example, in the solution to Problem I.1 he finds two numbers summing up to $100$ and having difference $40$; however, the method that he delineates always generalizes to arbitrary givens.

A present-day mathematician, when faced with an equation, would expect the author to specify the domain in which the solution should be sought. Diophantus on his part is not very explicit about this. He only makes two remarks bearing on his general concept of number (\textgreek{ἀριθμὸς} in Greek -- hence the title of the work) in the introduction to the {\it Arithmetica}. Right away in line 14 of the introduction, he states (\cite[p.~2]{tannery})
\begin{quotation}
\upshape
Now, you know that all numbers are composed {\it of some multitude of units} (\textgreek{ἐκ μονάδων πλήθους τινός})
\end{quotation} 
Later, after having defined squares, cubes, and so on (\cite[p. 6]{tannery}), he writes
\begin{quotation}
\upshape
But that which possesses none of these special properties, but contains in itself an {\it indeterminate multitude of units} (\textgreek{πλῆθος μονάδων ἀόριστον}), is called  a ``number'', and its sign is $\varsigma$.
\end{quotation} 

The ``sign'' mentioned is the sign for the unknown, as we will explain more fully in Section \ref{ss_notation}. The above remarks do not seem to give the modern reader much information about Diophantus' definition of numbers. Upon reading the {\it Arithmetica} however, it becomes clear that the solutions that Diophantus seeks for his problems are always what we would call {\it positive rational} numbers. Even in Book I, where the ``given'' quantities are usually such that the solutions consist of integers, it is clear that this is done just to keep the calculations simple (as Diophantus writes explicitly in the solution to I.24, for example).

Finally, it is worth mentioning the appearance of irrational quantities in the {\it Arithmetica}. Whether or not Diophantus considers these to be ``numbers'' (a question which Jacob Klein considers in his monumental study \cite{klein}, and which he answers with a decided negative) is not entirely clear, but occasionally he certainly recognizes irrational solutions to his equations. For example, in Problem IV.9 he encounters the equation $35x^2=5$, whereupon he remarks ``and $x$ becomes irrational'' (\textgreek{καί γίνεται ὁ $\varsigma$ οὐ ῥητός}) and then dismisses the solution. The safest conclusion seems to be that, whatever Diophantus' views were on the ontological status of irrational quantities, he certainly wanted the solutions to his {\it problems} to be rational.

\subsection{Algebraic notation in Diophantus' work} 
\label{ss_notation}
Throughout the {\it Arithmetica}, Diophantus employed a rudimentary form of algebraic notation. We will give an overview of his notation, not because we will use it, but in order to make the point that it enabled Diophantus to work with {\it univariate} polynomials in much the same way that we ourselves do.

Diophantus denotes the unknown by the symbol $\varsigma$; he never needs more than one at a time. The symbols for the powers of the unknown (up to the sixth) are given in the following table. Each symbol, with the possible exception of $\varsigma$, merely combines the first two letters of the corresponding Greek term, hence they can be regarded as abbrevations.\footnote{The symbol $\varsigma$ is only an approximate rendering of the one found in the manuscripts; for a discussion about its shape and provenance, see \cite{heath}.}

\begin{center}
\begin{tabular}{|c|c|c|c|}
\hline
symbol & name in Greek & literal meaning & our notation \\
\hline
\hline
$\stackrel{\textnormal{o}}{\operatorname{M}}$ & \textgreek{μονάς} & unit & $x^0$\\
\hline
$\varsigma$ & \textgreek{ἀριθμὸς} & number & $x^1$ \\
\hline
\vphantom{$\stackrel{\textnormal{o}}{\operatorname{M}}$} $\Delta^{\Upsilon}$ & \textgreek{δύναμις} & power, strength & $x^2$ \\
\hline
\vphantom{$\stackrel{\textnormal{o}}{\operatorname{M}}$}  $\textnormal{K}^{\Upsilon}$ & \textgreek{κύβος} & cube, die & $x^3$ \\
\hline
\vphantom{$\stackrel{\textnormal{o}}{\operatorname{M}}$}  $\Delta^{\Upsilon}\Delta$ & \textgreek{δυναμοδύναμις} & ``power-power'' &  $x^4$ \\
\hline
\vphantom{$\stackrel{\textnormal{o}}{\operatorname{M}}$}  $\Delta \textnormal{K}^{\Upsilon}$ & \textgreek{δυναμόκυβος} & ``power-cube'' & $x^5$ \\
\hline
\vphantom{$\stackrel{\textnormal{o}}{\operatorname{M}}$}  $\textnormal{K}^{\Upsilon}\textnormal{K}$ & \textgreek{κυβόκυβος} & ``cube-cube'' & $x^6$ \\
\hline
\end{tabular}
\end{center}

Here, the term \textgreek{δυναμοδύναμις}, or ``power-power'', is not to be read as ``the power of the power'', but rather as ``a power's worth of powers'', as is clear from the names for the fifth power (``a cube's worth of powers'') and the sixth power (``a cube's worth of cubes''). 

We mention here that the word \textgreek{τετράγωνος} ``square'' is frequently abbreviated by the symbol $\square$, but this should not really be regarded as a form of algebraic notation. The symbol can always be taken to mean the word ``square'' as if it were written out in full, and is never used to mean ``the square of the unknown''. (An example is given towards the end of this subsection.)

The Greek numerals are written using the letters of the Greek alphabet, but with a line drawn above them to indicate that they are to be interpreted as numbers. As a sample, we list the integers from $1$ to $10$, $100$ and $200$.

\begin{center}
\begin{tabular}{|c|c|}
\hline
Greek numeral & our notation \\
\hline
\hline
$\overline{\textnormal{\textgreek{α}}}$ &   $1$ \\
\hline
$\overline{\textnormal{\textgreek{β}}}$ &   $2$ \\
\hline
$\overline{\textnormal{\textgreek{γ}}}$ &   $3$ \\
\hline
$\overline{\textnormal{\textgreek{δ}}}$ &   $4$ \\
\hline
$\overline{\textnormal{\textgreek{ε}}}$ &   $5$ \\
\hline
$\overline{\textnormal{\textgreek{ς}}}$ &   $6$ \\
\hline
$\overline{\textnormal{\textgreek{ζ}}}$ &   $7$ \\
\hline
$\overline{\textnormal{\textgreek{η}}}$ &   $8$ \\
\hline
$\overline{\textnormal{\textgreek{θ}}}$ &   $9$ \\
\hline
$\overline{\textnormal{\textgreek{ι}}}$ &   $10$ \\
\hline
$\overline{\textnormal{\textgreek{ρ}}}$ &   $100$ \\
\hline
$\overline{\textnormal{\textgreek{σ}}}$ &   $200$ \\
\hline
\end{tabular}
\end{center}
Numbers which in decimal notation need more than one non-zero digit to write them down are written by simply concatenating the Greek numerals. So $208$ is written $\overline{\textnormal{\textgreek{ση}}}$, and $216$ is written $\overline{\textnormal{\textgreek{σις}}}$. Multiplicative inverses of both numbers and powers of the unknown are denoted by a symbol resembling $\chi$ in superscript, so that $\overline{\textnormal{\textgreek{δ}}}^\chi$ denotes $\dfrac{1}{4}$ and $\Delta^{\Upsilon \chi}$ denotes $\dfrac{1}{x^2}$.

The coefficient of each power of the unknown is written {\it after} the symbol for the power of the unknown; for example, $10 x^2$ would be written by Diophantus as $\Delta^\Upsilon \, \overline{\textnormal{\textgreek{ι}}}$. The symbol $\stackrel{\textnormal{o}}{\operatorname{M}}$, a ligature containing the first two letters of \textgreek{μονάς} (meaning ``unit'' or ``unity'') is used to signify the constant term of a polynomial, or simply to denote the number $1$; for example, Diophantus could write $\stackrel{\textnormal{o}}{\operatorname{M}} \overline{\textnormal{\textgreek{β}}}$ to mean simply the number $2$.

Addition of terms is denoted simply by juxtaposing them; so for example, $x+3$ is written as $\varsigma\, \overline{\textnormal{\textgreek{α}}} \stackrel{\textnormal{o}}{\operatorname{M}} \overline{\textnormal{\textgreek{γ}}}$. Subtraction of terms is denoted by a ``truncated inverted $\Psi$'', but for the purposes of this short overview we render it as simply an inverted $\Psi$ (\rotatebox[origin=c]{180}{$\Psi$}); moreover, since Diophantus does not have a separate sign for addition, the negative terms are written after the positive terms. So for example, the polynomial $9x^4-4x^3+6x^2-12x+1$ appearing in Problem IV.28 is written as
$$
\Delta^\Upsilon \Delta \, \overline{\textnormal{\textgreek{θ}}} ~
\Delta^\Upsilon \, \overline{\textnormal{\textgreek{ς}}} 
\stackrel{\textnormal{o}}{\operatorname{M}}  \overline{\textnormal{\textgreek{α}}} ~
\rotatebox[origin=c]{180}{$\Psi$} ~
\textnormal{K}^\Upsilon \, \overline{\textnormal{\textgreek{δ}}} ~
\mathbf{\varsigma} \, \overline{\textnormal{\textgreek{ιβ}}}.
$$

Finally, the abbrevation \textgreek{ἴσ}\textnormal{.} is sometimes used to denote equality. It is simply shorthand for \textgreek{ἴσος} ``equal'', or rather one of its declined forms, such as \textgreek{ἴσοι} or \textgreek{ἴσαι}. For example, the phrase ``$8x+4$ is equal to a square'' is written as
$$
\mathbf{\varsigma} \, \overline{\textnormal{\textgreek{η}}}  \stackrel{\textnormal{o}}{\operatorname{M}} \overline{\textnormal{\textgreek{δ}}} ~~ \textnormal{\textgreek{ἴσ}.} ~~ \square^{\textnormal{\textgreek{ῳ}}}.
$$

As explained, we will not use this notation in the rest of the article. However, it is hopefully clear at this point that using modern algebraic symbolism in translations of Diophantus' work, {\it when restricted to the use of one variable}, does not distort the original text in any essential way. Indeed, the overview given above should in principle enable the reader to reconstruct, from a univariate polynomial given in modern notation, the way in which Diophantus would have rendered it.

\subsection{On Diophantus and generality}
\label{dioph_gen}


As alluded to in the quotation by Euler, although Diophantus was not interested in giving general solutions to his equations, his methods can almost always be easily extended to give more general solutions, sometimes even {\it the} general solution. As an example of this, I will give a translation of the very straightforward Problem I.14 (which is the first ``indeterminate'' problem of the {\it Arithmetica}\footnote{This nomenclature is traditional in the literature on Diophantus, cf. \cite{hankel}. The distinction ``determinate'' versus ``indeterminate'' can be taken to mean that the algebraic variety associated to the equations has dimension zero or greater than zero respectively. Moreover, in almost all cases in the {\it Arithmetica}, the set of rational solutions of an ``indeterminate'' equation is infinite, whereas that of a ``determinate'' equation is always finite.}) and its solution. 
\begin{quotation}
{\sc Problem I.14.} To find two numbers such that their product has to their sum a given ratio.
\end{quotation}
We now give a literal translation of Diophantus' solution. He starts off by stating a condition that has to be fulfilled by the solution, hereafter called a {\it diorism}\footnote{Often, the diorism imposes conditions on the ``given'' numbers, or ratios, mentioned in the problem, whereas in the present case a property of the {\it solution} is described. Diophantus uses the word \textgreek{προσδιορισμός}, but we prefer the term {\it diorism}, which is more common in ancient Greek mathematics, and whose meaning does not appear to be much different.}:
\begin{quotation}
We will assume a value for one of the numbers; it must be the case that this assumed value is greater than the given ratio.\footnote{This is to ensure that the resulting solution will be in positive numbers.} 

Let it be prescribed, then, that the product has to the sum the ratio $3$. 

Let one of the unknowns be $x$; the other, according to the diorism, must be greater than $3$. Let it be $12$. Then their product is $12x$, and their sum $x+12$. What remains is for $12x$ to be equal to three times $x+12$, since three times the smaller is equal to the greater. And $x$ becomes $4$.

One of them shall be $4$, the other $12$. And they solve the problem.
\end{quotation}

There are two features of Diophantus' solution that will strike the modern reader as peculiar: firstly, while Diophantus asks for a solution in the case of an arbitrarily ``given'' ratio, he only gives a solution in the case where the ratio is $3:1$. Secondly, he assumes the value $12$ for one of the unknowns, seemingly at random. 

The first feature is typical of the {\it Arithmetica}. Clearly, Diophantus' notational apparatus is ill-suited for solving arithmetic problems in their proper generality, and one can imagine that if Diophantus would have been forced to give a properly general treatment, such as would be required of a present-day mathematician, his presentation would have become much less digestible. But aside from this, it is clear from many passages in the {\it Arithmetica} that Diophantus' text must be construed as ``instruction by example''. The solutions in the {\it Arithmetica} are to be taken {\it mutatis mutandis}, as when one draws a specific triangle in the course of proving a fact true for {\it all} triangles.\footnote{\label{entoiaoristoi}Diophantus does on occasion explicitly ask for, and describe, one-parameter families of solutions (so for example in Problem IV.19), but this invariably comes at a considerable cost in readability. See also the solution to Problem V.29 (Section \ref{secV29}), where at one point Diophantus is effectively dealing with three unknowns at a time (which I have there labeled $p$, $q$, and $m$ for the sake of clarity), but expresses his argument entirely in words.}

The other remarkable feature is the sudden choice of $12$ for the second unknown. With this choice, the problem becomes ``determinate'', that is, it now admits of only a single solution. It is however clear to anybody who takes the trouble of replicating Diophantus' argument with a different number instead of $12$, that any other positive rational number greater than $3$ would serve equally well. Diophantus so much as indicates this by prefixing the choice by the consideration that the unknown has to be chosen larger than $3$. Other than that, seems to be the implication, there are no restrictions. (Here, one is somehow reminded of the ancient Roman legal principle of {\it exceptio probat regulam}, ``the exception proves the rule,'' which meant: wherever an exception is mentioned, one may surmise the existence of a general rule to which the exception applies. In our case, Diophantus excepts values of $3$ or less for the second unknown, so one naturally assumes that all other values are allowed.) As before, we can conclude that, while there is a nominal loss in generality, the generality is nevertheless clearly implied by the text.

\subsection{``Just a single solution''}

In almost\footnote{That Diophantus understood that some of his equations had infinitely many solutions is clearly shown in his treatments of Problems III.12 and III.19. In the latter, he states that ``we have seen how to divide a given square into two squares in infinitely many ways'' (\textgreek{ἐμάθομεν τόν δοθέντα $\square^{\textnormal{\textgreek{ον}}}$ διελεῖν εἰς δύο $\square^{\textnormal{\textgreek{ους}}}$ ἀπειραχῶς}), and indeed, this is a corollary to his solution of Problem II.8. However, remarks such as these are never made for their own benefit, but because of their relevance to the problem at hand: in the solution to III.19, for instance, he needs to write $65^2$ as a sum of two squares in not one but four different ways, which makes the infinitude of solutions to Problem II.8 a pertinent fact. For more on this, see \cite{schappacher}.} all problems of the {\it Arithmetica}, Diophantus limits himself to giving just {\it one} solution to his problem -- we saw as much above in the discussion of Problem I.14. However, as we argued there, the methods by which Diophantus arrives at his solutions can in most cases be adapted to give infinitely many solutions. As the great Leonhard Euler writes\footnote{The translation is from Heath \cite{heath}.}

\begin{quotation}
Diophantus himself, it is true, gives only the most special solutions of all the questions which he treats, and he is generally content with indicating numbers which furnish one single solution. But it must not be supposed that his method was restricted to those very special solutions. In his time the use of letters to denote undetermined numbers was not yet established, and consequently the more general solutions which we are now enabled to give by means of such notation could not be expected from him. Nevertheless, the actual methods which he uses for solving any of his problems are as general as those which are in use today; nay, we are obliged to admit that there is hardly any method yet invented in this kind of analysis of which there are not sufficiently distinct traces to be discovered in Diophantus. 
\end{quotation}

However, the fact remains that Diophantus, in the great majority of his solutions, though not in all, was satisfied with giving a single solution to his problems. This fact is the basis for the following judgment from Hankel's influential work \cite{hankel}, which is echoed in many later sources. Hankel writes \cite[p. 166]{hankel}

\begin{quotation}
But that which most of all robs his -- often surprisingly clever -- solutions of their scientific worth, is the fact that Diophantus invariably contents himself with selecting just a single solution to any one of his problems, out of the great mass of solutions that it might have,  without so much as glancing at any other.
\end{quotation}
It is one of the purposes of this article to show just how misguided this judgment is, even if we disregard the instances where Diophantus, contrary to Hankel's claim, actually {\it does} give one-parameter families of solutions explicitly (see footnote \ref{entoiaoristoi}). The fact is that Diophantus' methods can often easily be varied to give additional solutions -- we saw an example of this in Section \ref{dioph_gen}, and we will see many more -- and this alone means that he {\it does} glance at other solutions! He does not bother to spell out the details, but it can be done, and -- to the mind of the author -- it {\it should} be done if one wants to gain a realistic assessment of Diophantus' work. 

\subsection{Diophantus' solutions in terms of modern algebra}

In the present article, I have tried to make the ``hidden generality'' in Diophantus' work visible in six special cases. I present six problems from Diophantus' {\it Arithmetica}, together with either translations or synopses of his original solutions, stated in a simple algebraic language that is close to that of Diophantus himself. The original solutions are then followed by their ``modern translations'', where the full force of modern algebraic symbolism is brought to bear.

In this latter process of translation, I have strictly adhered to the following two rules:
\begin{itemize}
\item[(i)] at each point where Diophantus makes a seemingly completely arbitrary choice, substitute a ``general parameter'' instead of a specific constant;
\item[(ii)] at each point where Diophantus chooses a fixed quantity {\it subject to a condition alluded to in the original text}, I consider: first, whether that condition results in an arithmetical problem that Diophantus himself could have solved, with the evidence being the contents of the {\it Arithmetica} itself; second, in which generality he could have solved it; and finally, third, proceed by substituting the corresponding ``generic solution'' to the problem mentioned.
\end{itemize} 
Note that an application of rule (ii) starts a recursion with rule (i) as its base case, since that is the only case where the actual substitution of
 parameters occurs.

To summarize: it is my view that Hankel's judgment on that Diophantus ``invariably contents himself with [...] just a single solution [...] without so much glancing at any other'' does Diophantus an injustice. This injustice is the greater for the influence that Hankel's judgment has had on later commentators. I think that, by explicating Diophantus' solutions according to the two rules given above, we will be able to see the true scope and strength of Diophantus' work. 

\subsection{Sources and acknowledgments}

My acquaintance with the {\it Ar\-ithm\-et\-ic\-a} was first made through the excerpts given in the two-volume anthology by Ivor Thomas \cite{thomas}, which has appeared as part of the Loeb Classical Library (Volumes 335 and 362), and which contains the Greek originals and English translations side by side. I have also consulted the classical English edition of the {\it Arithmetica} by Sir Thomas L. Heath \cite{heath}, which includes a very thorough study of both the form and the content of the work, as well as a wealth of material drawn from other commentaries.\footnote{It must be said here that Heath's work is perhaps less a translation of Diophantus than a synopsis. Heath does an extremely good job in summarizing the mathematical ideas present in each solution, but he does leave the reader wondering in what exact way Diophantus himself phrased his solutions.} My translations given in this article, though they are inevitably influenced by the works of Thomas and Heath, are all done from what is still considered the ``official'' edition by Paul Tannery \cite{tannery}, which contains the original Greek with a Latin translation. 

I thank Oliver Braunling, Michael Gr\"ochenig, Francesco Polizzi, Roy Schepers, Frits Veerman, Allard Veldman, and Jesse Wolfson for their assistance, encouragement, and useful discussions.

\section{Problem II.20: a del Pezzo surface of degree $4$}

The twentieth problem of Book II runs thus:
\begin{quotation}
To find two numbers such that the square of either when added to the other yields a square.
\end{quotation}

\subsection{Diophantus' solution to Problem II.20}
The following is a translation of Diophantus' solution, writing $x$ for his unknown. 
\begin{quotation}
Let the first number be fixed as $x$, and the second as $2x+1$, in order that the first requirement be satisfied. The second requirement then is that
$$
x+(2x+1)^2 = 4x^2+5x+1
$$
be square. Let it be the square of $2x-2$. Then we have
$$
4x^2+5x+1 = (2x-2)^2 = 4x^2-8x+4 ~ \Longrightarrow ~ 13x = 3,
$$
which yields $x = \dfrac{3}{13}$. Hence the first number is $\dfrac{3}{13}$, and the second $\dfrac{19}{13}$.
\end{quotation}

\subsection{Comments on the solution to Problem II.20}

In modern notation, Diophantus asks for (positive) rational numbers $x,y$, such that there exist further rational numbers $u,v$ such that
\begin{equation}
\label{II20}
\left\{ ~~ \begin{array}{lcl} x^2 + y & = & u^2 \\ x + y^2 & = & v^2 \end{array} \right.
\end{equation}

The present-day reader can hardly fail to notice that Diophantus makes two choices in his solution that seem rather arbitrary to one familiar with modern algebra. To begin with, in his first step, he assumes that the indeterminate $y$ equals some linear expression in $x$, to wit $y = 2x+1$, which has the virtue of making $x^2+y$ the square of a linear polynomial in $x$. Of course, any other expression of the form $y = 2\lambda x + \lambda^2$ would have served him equally well here. Similarly, when he takes $v = 2x-2$, he could have chosen $v = 2x + \mu$ for any $\mu$; indeed, it would have seemed much less strange to modern eyes if he had. 

We will now reinterpret Diophantus' solution in terms of algebraic geometry. We will show that in doing so, we end up with a birational parametrization of the algebraic surface defined by \eqref{II20}. Thus we will see that, by recasting Diophantus' solution to Problem II.20 in terms of modern algebra, one ends up with a far more general solution to the problem than Diophantus did himself. 

\subsection{The geometry behind the solution to Problem II.20}

The equations \eqref{II20} define a smooth intersection of two quadrics in $\mathbb{P}^4$, hence a del Pezzo surface $X$ of degree 4 over the field $\mathbb{Q}$ of rational numbers \cite[Proposition IV.16]{beauville}. We replace Diophantus' substitution $y=2x+1$ with the more general $y=2\lambda x + \lambda^2$ for some indeterminate $\lambda$ (with $\lambda=1$ corresponding to the choice made by Diophantus). Substituting this into the second equation gives
$$
x + (2\lambda x + \lambda^2)^2 = 4 \lambda^2 x^2 + (4 \lambda^3 + 1)x + \lambda^4 = v^2.
$$
At this point, Diophantus' substitution $v=2x-2$ (which works for his value of $\lambda$) was designed to make the square terms on both sides cancel. To accomplish this, the modern mathematician would rather choose $v = 2\lambda x + \mu$, with $\mu$ again indeterminate, since this accomplishes the same cancellation of terms but allows for a more general solution. The resulting equation for $x$ is
$$
4 \lambda^2 x^2 + (4 \lambda^3 + 1)x + \lambda^4 = (2 \lambda x + \mu)^2,
$$
giving
$$
(4 \lambda^3 + 1)x+\lambda^4 = 4\lambda\mu x+\mu^2 ~ \Longrightarrow ~ x = \frac{\mu^2-\lambda^4}{4 \lambda^3 - 4\lambda\mu  + 1}.
$$
Together with $y=2\lambda x + \lambda^2$, $u = x + \lambda$, and $v = 2\lambda x + \mu$, this defines a rational map
$$
\phi \colon \mathbb{P}^2 \dashrightarrow X,
$$
where $\lambda,\mu$ are the coordinates on  $\mathbb{P}^2$. The rational map $\phi$ is even birational, its inverse being given by $\lambda = u - x$, $\mu = v - 2\lambda x$.

To summarize: we started from Diophantus' solution to \eqref{II20}, we replaced two seemingly arbitrary choices of specific rational numbers by analogous substitutions of the indeterminates $\lambda$ and $\mu$, and ended up with a birational parametrization of the del Pezzo surface of degree 4 defined by \eqref{II20}. 

\subsection{A conic bundle structure on $X$} 
We can give a further geometric interpretation to Diophantus' manipulations.  The substitution $y = 2\lambda x + \lambda^2$ into the second equation \eqref{II20} means, in more geometric terms, that we restrict to looking at the (reducible) curve $u = \pm(x + \lambda)$ on $X$. Conversely, by choosing a sign in the latter expression, say 
\begin{equation}
\label{II20magic}
u=x+\lambda,
\end{equation}
it follows from the first equation that we have $y = 2\lambda x + \lambda^2$. Applying the substitution \eqref{II20magic} for an indeterminate $\lambda$ effectively means endowing $X$ with the structure of a conic fibration, as follows
\begin{align*}
\pi  \colon X & \dashrightarrow \mathbb{P}^1 \\\
(x,y,u,v) & \longmapsto u-x
\end{align*}
To describe the first step of Diophantus' solution in these terms actually helps in understanding the second step as well, since the remainder of Diophantus' argument is essentially to find a rational point on the fibre of $\pi$ over $\lambda = 1$. In our more geometric rendering of his argument, we considered the {\it generic} fibre of $\pi$, which is a conic $C$ over the field $\mathbb{Q}(\lambda)$, and this likewise has a rational point. This just rephrases the fact that $\pi$ has a section, so that it is once more clear that $X$ is birationally equivalent to $\mathbb{P}^2$.

\section{Problem II.31: a singular rational surface}
The thirty-first problem of Book II is as follows:

\begin{quotation}
To find two numbers whose sum is a square, such that their product plus or minus their sum gives a square.
\end{quotation}

\subsection{Diophantus' solution of Problem II.31}
The following is a translation of Diophantus' solution:
\begin{quotation}
We first observe that the numbers $2$ and $4$ have the property that twice their product, which is $16$, is a square; also, the sum of their squares, which is $20$, plus or minus twice their product is square.\footnote{\label{Bachet}The original text of the {\it Arithmetica} is somewhat confusing at this point; our translation is partly based on the emendation proposed by Bachet. I will here give a more faithful translation of Diophantus' text, and add an explanation. The first paragraph of the original solution reads thus: ``Because, if there are two numbers of which one is twice the other, the sum of their squares after either subtracting or adding twice their product gives a square, we put forward $4$ and $2$.'' That is, if $a$ and $b$ are such that $b=2a$, then both $a^2+b^2-2ab$ and $a^2+b^2+2ab$ are squares (and Diophantus proposes to take $a=2,b=4$). This is a strange claim, mainly because the condition $b=2a$ seems out of place, and also because it is perhaps not immediately clear how this claim connects up with the rest of the solution. The claim as it stands is rather trivial, since $a^2+b^2\pm 2ab = (a \pm b)^2$, which Euclid already knew, and Diophantus must have realized that the condition $b=2a$ is unnecessary. The mention of this condition probably points towards the strategy that becomes clear with his next move, namely to assume the product of the two unknowns to be equal to $a^2+b^2$ and the sum to $2ab$ (or, as it turns out, the product to $(a^2+b^2)x^2$ and the sum to $2ab x^2$, but this seems to be a further, independent move), so that the last two requirements of the problem would be met. But in doing so, he needs $2ab$ to be a square in order to meet the {\it first} requirement, and this is where the condition $b=2a$ comes in. So this seems to be Diophantus' strategy, and he takes $a=2$ and $b=4$, which gives $a^2+b^2=20$ and $2ab=16$, so that the product of the two unknowns is to become $20x^2$, and their sum $16x^2$. (In his 1621 edition of the {\it Arithmetica}, Bachet inserts an explanation to this effect in the original text, which Tannery later included in his edition as a footnote \cite{tannery}. The translation given in the main text partially follows Bachet's text.)} 

Therefore, multiplying throughout by $x^2$, let the product of the unknowns be $20x^2$, and their sum $16x^2$ [so that the requirements are met]. Let one be $2x$, and the other $10x$. Their sum is then $12x$, but also $16x^2$.

We then have $16x^2=12x$, hence $x=\dfrac{12}{16} = \dfrac{3}{4}$. Hence the first number is $2x = \dfrac{3}{2}$, and the second $10x = \dfrac{15}{2}$, which solve the problem.
\end{quotation}

The above solution, bizarre though it is, contains enough ideas to deduce from it easily a one-to-one parametrization of the associated surface. This we will explain next.

\subsection{Comments on the solution to Problem II.31}
\label{II31comments}

The solution given by Diophantus is extremely brief, and some amplification is in order. The actual idea of the solution, which is very typical for the {\it Arithmetica}, becomes clear in the second paragraph: Diophantus will impose values on the product $P$ and the sum $S$ of the two unknowns, which he will design in such a way that the requirements of the problem are met; that is, 
\begin{equation}
\label{II31_PandS}
\textnormal{$S$ is a square and $P \pm S$ are both squares.}
\end{equation}
However, this idea as it stands is not quite good enough, since a system of the form
\begin{equation}
\label{II31a_pre}
\left\{ ~~ \begin{array}{lcl} xy & = & P \\ x+y & = & S  \end{array} \right.
\end{equation}
is not always soluble in rational numbers.\footnote{In fact, Problem I.27 asks precisely for the solution of the system \eqref{II31a_pre}, and it is noted there that the rationality of $x$ and $y$ is equivalent to 
$$
\dfrac{1}{4} S^2 -P
$$
being a square.} Thus Diophantus uses a trick that he will use many times later in the work, namely that of playing off (weighted) homogeneous and inhomogeneous equations against each other. Since the conditions \eqref{II31_PandS} are such that $P$ and $S$ can be scaled by the square of an additional parameter without altering the truth value of the statements, Diophantus lets the square of his unknown act as scaling factor. Denoting the unknown here by $t$, this means that he ends up with a system of the form
\begin{equation}
\label{II31a_pre_gen}
\left\{ ~~ \begin{array}{lcl} xy & = & Pt^2 \\ x+y & = & St^2 \end{array} \right.
\end{equation}
which he solves by writing $P=QR$ for rational numbers $Q$ and $R$ (in his case $P=20$, and he takes $Q=2$ and $R=10$) and assuming $x=Qt$ and $y=Rt$. In this way, the first equation is satisfied, and the second yields
$$
Qt+Rt =St^2, 
$$
which leaves him with a linear equation for the unknown $t$.

As for the subproblem of finding suitable values $P$ and $S$, the text is somewhat elliptical (as indicated in footnote \ref{Bachet}), but the intended strategy is  quite clear. There are two main ideas. The first is to use the identities
$$
a^2+b^2\pm 2ab = (a \pm b)^2,
$$
and the second is to choose $a$ and $b$ such that $b=2a$. Now the  latter relation ensures that  $2ab$ will be a square, so that  we may put $P=a^2+b^2$ and $S=2ab$. Again, it is striking to see how Diophantus sets up the solution to his subproblem in such a way that it is almost impossible not to notice how it might be generalized: if $a$ is given, then for any rational number $r$ one may set $b=2r^2a$, and it follows that $2ab=(2ra)^2$ is square. Hence, for any $a$ and $r$, one may choose
\begin{equation}
\label{II31forlater}
P =a^2+b^2 = (4r^4+1)a^2,~~ S=4r^2a^2 
\end{equation}
for substitution into \eqref{II31a_pre_gen}.

\subsection{Geometrical interpretation of the solution}
Problem II.31 leads to a surface $X$ given by 
\begin{equation}
\label{II31a}
\left\{ ~~ \begin{array}{lcl} x + y & = & u^2 \\ xy + x + y & = & v^2 \\ xy - x - y & = & w^2 \end{array} \right.
\end{equation}

This is a complete intersection of three quadrics in $\mathbb{P}^5$, with isolated non-ordinary singularities (as calculated with \texttt{magma} \cite{magma}). We will see later that $X$ is in fact a rational surface. In fact, following Diophantus' solutions, we will see that $X$ can be fibred in (singular) rational curves, with a conic having a rational point as the base curve; moreover, the fibration has a section. In sum, Diophantus' solution prepares the way for a birational parametrization of $X$.

\subsubsection{A rational curve on $X$}
The ``Ansatz'' proposed by Diophantus is to assume a value for the product $xy$, namely $20t^2$, and one for the sum $x+y$, namely $16t^2$. This yields the curve
\begin{equation}
\label{II31fibre}
C \colon xy = 20t^2,~~x+y = 16t^2
\end{equation}
that is an intersection of two quadrics in $\mathbb{P}^3$ with a singularity at $(0,0,0)$, so that $C$ is rational. For the remaining variables we may take, for example,
\begin{equation}
\label{II31fibremore}
u = 4t , ~~ v = 6t, ~~ w = 2t,
\end{equation}
other sign choices being also possible. Taken together, \eqref{II31fibre} and \eqref{II31fibremore} define a singular rational curve on $X$. 
Diophantus solves the system \eqref{II31fibre} by setting
$$
x = 2t, ~~ y = 10t.
$$
There is clearly nothing special about the choice of coefficients $2$ and $10$ other than their product being $20$, so we could just as well take
\begin{equation}\label{II31pencil}
x = \lambda t, ~~ y = \dfrac{20t}{\lambda}.
\end{equation}
In geometrical terms, \eqref{II31pencil} defines a one-dimensional pencil $\{L_\lambda\}_\lambda$ of lines on the quadric surface
$$
Q \colon xy = 20t^2
$$
defined by the first equation of the system \eqref{II31fibre}. Note that all of the lines of this pencil pass through the singular point $(0,0,0)$ of $Q$. A general line $L_\lambda$ of the pencil will intersect the quadric defined by the second equation in two points. Such a $L_\lambda$ will therefore intersect $C$ in a unique rational point $P$ besides $(0,0,0)$. Going through the calculations, we obtain for the coordinates of $P$ the following values
$$
x = \dfrac{20 +\lambda^2}{ 16 },~~ y = \dfrac{100 +5\lambda^2}{ 4 \lambda^2}, ~~t = \dfrac{20 +\lambda^2}{ 16 \lambda}.
$$ 

\subsubsection{A pencil of rational curves on $X$} So far then, we see that Diophantus' ideas have yielded the par\-am\-et\-riz\-at\-ion of a rational curve on $X$. However, the choice of the coefficients $20$ and $16$ appearing in \eqref{II31fibre} still seems arbitrary. Indeed, despite its obscurities, the text of Problem II.31 gives strong hints towards a generalization. As we saw in the discussion leading up to \eqref{II31_PandS}, the choice of the coefficients $20$ and $16$ arose as a special case of a more general problem (see footnote \ref{Bachet}), namely that of finding rational numbers $P$ and $S$ such that 
\begin{equation}
\label{II31basic_problem}
S = L^2, ~~ P + S = M^2, ~~ P - S = N^2
\end{equation}
for certain rational numbers $L$, $M$, and $N$. 
Eliminating $P$ and $S$, we find that the problem is equivalent to finding the rational points on the projective conic
$$
B \colon M^2-N^2=2L^2.
$$
Writing $m=M/L$ and $n=N/L$, we get the affine equation  
$$
B \colon m^2-n^2=2.
$$
Each point $\beta = (m_0,n_0)$ on $B$ corresponds to a curve $X_\beta$ on $X$ given by the equations
\begin{equation}
\label{II31superfibre}
\left\{ 
\begin{array}{lcl}
x+y & = & t^2  \\ 
xy+x+y & = & m_0^2 t^2  \\ 
xy-x-y & = & n_0^2 t^2 \\
u & = & t \\
v & = & m_0 t \\
w & = & n_0 t
\end{array}
\right.
\end{equation}
where of course the $t$ is auxiliary (that is, the curve $X_\beta$ is obtained by projecting away from the $t$-coordinate). We see at once from \eqref{II31superfibre} that, for any point $(x_0:y_0:u_0:v_0:w_0)$ on $X$, the corresponding point $\beta \in B$ can be recovered by $m_0 = v_0/u_0$ and $n_0 = w_0/u_0$.

In sum, we can define a rational map from $X$ to $B$, as follows
\begin{align*}
\pi \colon X & \dashrightarrow B \\\
(x,y,u,v,w) & \longmapsto (v/u,w/u)
\end{align*}
which gives $X$ the structure of a fibration into rational curves, with the rational curve $B$ as basis. (Incidentally, this shows that $X$ is a rational surface, as claimed before.) Moreover, the original argument given by Diophantus can be replicated to show that $\pi$ even comes with a section. Considering $X$ as a rational curve over the function field $K$ of $B$, and imitating the construction from \eqref{II31fibre} onwards, we project the generic fibre $X_K$ down to the following intersection of two quadrics in $(x,y,t)$-space
$$
C \colon \, 
\left\{
\begin{array}{lcl}
xy &  =  & (n^2+1) t^2\\
x + y & = & t^2
\end{array}
\right.
$$
We consider the following pencil of lines
$$
L_\lambda \colon ~ x = \lambda t, ~~ y = \dfrac{(n^2+1) t}{\lambda}.
$$
Intersecting $L_\lambda$ with $C$ yields 
\begin{equation}\label{II31preparam}
x = \lambda^2+n^2+1, ~~y = n^2+1 +\left( \dfrac{n^2+1}{\lambda}\right)^2 , ~~t= \dfrac{\lambda^2 + n^2+1}{\lambda}.
\end{equation}
Finally, if we would also parametrize $B$, we can use these formulas to give a parametrization of $X$. We will do this in the final subsection, again making crucial use of ideas in Diophantus' text.

\subsubsection{Rational points on $B$}
To recapitulate, we have seen how Problem II.31 leads to a rational surface $X$, and how Diophantus' solution naturally leads one to consider $X$ as a fibration $\pi \colon X \rightarrow B$ into rational curves, where $B$ is the rational curve given by
$$
B \colon m^2-n^2=2,
$$
and that it even leads to the construction of a section of $\pi$. 

We further recall that the equation for $B$ arose from the problem of finding $P$ and $S$ such that $S$, $S+P$, and $S-P$ are all squares of rational numbers. We work with an affine equation for $B$, which is equivalent to scaling $P$ and $S$ so that $S=1$; we then also have $P+S=m^2$ and $P-S=n^2$. 

To the modern eye, it is of course a complete triviality that $B$ has rational points, since the left-hand side of its equation factors as
$$
m^2-n^2=(m+n)(m-n),
$$
which shows that we may parametrize $B$ by putting, for example, $m+n=2a$, $m-n=1/a$, and solving for $m$ and $n$ to get $(m,n)=(a+1/(2a),a-1/(2a))$. However, in keeping with the spirit of this article, we would prefer to consider  the idea that Diophantus himself proposes. His text clearly alludes to the identities
\begin{equation}\label{II31anothergreatidea}
a^2+b^2\pm 2ab = (a \pm b)^2,
\end{equation}
and to the idea of setting $S=2ab$, $S+P=(a+b)^2$, and $S-P=(a-b)^2$. Since we have scaled things so that $S=1$, we have to add the requirement that $2ab=1$, so that $b=1/(2a)$. With these values, equation \eqref{II31anothergreatidea} leads to $S=1$ and $P=\dfrac{4a^4+1}{4a^2}$, which of course is the special case of \eqref{II31forlater} where $r=1/(2a)$, and to 
$$
S = 1,~~P + S = \dfrac{4a^4+4a^2+1}{4a^2} = \left( a+\dfrac{1}{2a} \right)^2=m^2,
$$
and
$$
P-S = \dfrac{4a^4-4a^2+1}{4a^2} = \left( a-\dfrac{1}{2a} \right)^2=n^2.
$$
This immediately yields the following parametrization of $B$:
$$
(m,n) = \left( a+\dfrac{1}{2a}, a-\dfrac{1}{2a} \right),
$$
which is the same parametrization as the one found above. Note that the parametrization is birational, since $m+n=2a$, so that $a$ can be recovered from $(m,n)$. 

Combined with \eqref{II31preparam}, this yields the following birational parametrization of $X$:
\begin{equation*}
x = \lambda^2 + a^2 + \dfrac{1}{4a^2}, ~~  y = \left( a^2+\dfrac{1}{4a^2} \right)  \dfrac{x}{\lambda^2}, ~~ u = \dfrac{x}{\lambda},
\end{equation*}
and
\begin{equation*}
v = \left( a+\dfrac{1}{2a} \right) \dfrac{x}{\lambda} , ~~  w = \left( a-\dfrac{1}{2a} \right) \dfrac{x}{\lambda} .
\end{equation*}

\section{On intersections of two quadrics in $\mathbb{P}^3$ in the ``Arithmetica''}
\label{double_eqs}

In many problems, Diophantus is confronted with what he calls a {\it double equation}\footnote{The terms used by Diophantus are \textgreek{διπλοϊσότης}, \textgreek{διπλῆ ἰσότης} (both occurring e.g. in II.11), and \textgreek{διπλῆ ἴσωσις} (e.g. in II.13).}. We will encounter an example of such a double equation in the next section. Stated in modern language, this type of problem comes down to determining a rational point on a curve of genus zero or one, given as an intersection of two quadrics in $\mathbb{P}^3$. In this section, we give an overview of the double equations featured in the {\it Arithmetica}, and we discuss Diophantus' methods for dealing with the genus one case.\footnote{My main authority for the overview given in this section besides the {\it Arithmetica} itself is the outstanding discussion by Heath \cite[pp. 81--87]{heath}. Heath gives a concise but complete classification of the so-called ``double-equations of the second order'' that appear in the {\it Arithmetica}. One difference between our presentation and that of Heath is that that of us is more geometrical, while his is purely algebraic. Unsurprisingly, our geometric approach sheds some further light on the subdivision found in Heath's work. Heath distinguished three cases, roughly according to (i) $a_1 \neq 0$ and $a_2 \neq 0$; (ii) $a_1 \neq 0$ and $a_2 = 0$; and (iii) $c_1=c_2=0$ (which after a projective transformation becomes equivalent to $a_1=a_2=0$). The distinction between (i) and (ii), here dealt with in Section \ref{case_dist}, is fairly subtle, even in Heath, but Diophantus' solution method for type (iii) is wholly different from the other cases. Geometrically, this makes perfect sense: equations of type (iii) lead to curves of (geometric) genus zero, while equations of type (i) and (ii) yield curves of genus one in general. Even so, Diophantus' over-all strategy does seem to derive from his solution for the genus zero case, as we shall see shortly.}

\subsection{Definition} Diophantus calls a {\it double equation} a system of two equations in which two polynomials in $x$ have to be made into either a square or a cube. That is, a double equation for Diophantus is a system of the form
\begin{equation}
\label{double_eq}
\left\{ ~~ \begin{array}{lcl} p_1(x) & = & u^m \\ p_2(x) & = & v^n \end{array} \right.
\end{equation}
where $p_1(x)$ and $p_2(x)$ are polynomials in $x$ with rational coefficients, and where $m$ and $n$ are integers (specified beforehand), equal to either $2$ or $3$, although usually both equal to $2$. In this section, we will restrict to $m=n=2$, which yield the equations referred to by Heath as ``double-equations of the second order'' \cite{heath}, or in other words, intersections of two quadrics in $\mathbb{P}^3$.

In the cases occurring in the {\it Arithmetica}, the polynomials $p_1$ and $p_2$ are always at most of degree $2$. That is, the system takes the following form:
\begin{equation}
\label{double_eq_expl}
\left\{ ~~ \begin{array}{lcl} a_1 x^2+b_1 x+c_1 & = & u^2 \\ a_2x^2+b_2x+c_2 & = & v^2 \end{array} \right.
\end{equation}
so that the associated projective curve $C$ is of genus zero or one.\footnote{\label{smoothness_crit}It is easily verified that $C$ is irreducible if and only if $p_1,p_2$ are non-constant and without repeated roots and $p_1/p_2$ is non-constant, and that $C$ is smooth if and only if the polynomial $p_1 p_2$ is of degree $3$ or $4$ and has no repeated roots.} Diophantus considers a wide variety of special cases of \eqref{double_eq_expl}, but he has a very uniform strategy for dealing with them that he follows in most cases.\footnote{The exceptional cases are the ones where either $a_1=a_2=0$, which is the most common kind of double equation appearing in the {\it Arithmetica}, or $c_1=c_2=0$. In both cases, $p_1$ and $p_2$ have a zero in common, hence $C$ is singular, and so it is rather satisfying from a modern point of view that Diophantus treats these cases differently from the general case. Heath \cite{heath} classifies the case where $a_1=a_2=0$ as a ``double-equation of the first order'', while treating the one where $c_1=c_2=0$ as being of the second order. The latter case occurs in Problems VI.12 and VI.14, where Diophantus reduces the system in short order to a plane conic \cite{heath}.} We give a description of his method in this section. 

\subsection{Method of solution} The method that Diophantus uses to solve his double equations is directly based on his treatment of the case where $a_1=a_2=0$, in which case $C$ is of genus zero. For this case he gives a clear method in his solution to Problem II.11 (in which the term ``double equation'' is used for the first time):
\begin{quotation}
\upshape
{\sc Problem II.11.} To add to two given numbers the same number, and make each into a square.
\end{quotation}
We give the part of the solution that is relevant to us:
\begin{quotation} 
\upshape
Let the given numbers be $2$ and $3$, and let $x$ be the number to be added. Then $x+2$ on the one hand, and $x+3$ on the other, will be squares. This form is called a {\it double equation}, and it is solved in this way: consider the difference, and seek two numbers whose product is this difference: they are $4$ and $\dfrac{1}{4}$. Of these two numbers, take half of their difference squared to be the smaller square, and half their sum squared to be the larger square. 

But the square of half their difference is $\dfrac{225}{64}$. This is to be equal to $x+2$, and $x$ becomes $\dfrac{97}{64}$. [...]
\end{quotation}

It is clear that the ``two numbers whose product is this difference'' may be chosen arbitrarily, which reflects the fact that the curve $C$ associated to the double equation is of genus zero. In the general case, however, the curve $C$ will be of genus one, and the two numbers have to be chosen in a special way, as we will see shortly.

We will now assume that the system \eqref{double_eq_expl} is properly ``of the second order'' (in the  terminology of Heath), that is, $a_1,a_2$ are not both zero. We will describe Diophantus' strategy. Along the way, we will encounter the conditions under which it works.\footnote{We will see that there are two major necessary conditions, being the existence of the factorization \eqref{crucial_factorization} and the existence of a suitable rational value for the parameter $\lambda$ occurring in \eqref{subst}.
} As in his solution to Problem II.11, he starts out by subtracting one equation from the other, obtaining 
\begin{equation}
\label{subtractsecondfromfirst}
u^2 - v^2 = (a_1-a_2) x^2 + (b_1 - b_2) x + (c_1 - c_2).
\end{equation}
It is crucially important for Diophantus' strategy to work that the right-hand side of this equation splits into linear factors over the rationals. (For an interpretation of this in geometric terms, see below.) Let us assume that this is the case, and put
\begin{equation}
\label{crucial_factorization}
(u+v)(u-v) = u^2 - v^2 = (m_1 x + n_1)(m_2x + n_2),
\end{equation}
where the $m_i$ and $n_i$ are rational numbers. Analogously to his solution of Problem II.11, Diophantus proceeds to equate the pairs of linear factors appearing in \eqref{crucial_factorization} to each other, up to a well-chosen constant $\lambda$. That is, he puts
\begin{equation}
\label{subst}
\left\{ ~~ \begin{array}{lcl} u+v & = & \lambda \cdot (m_1 x + n_1) \\ u-v & = & \lambda^{-1} \cdot (m_2 x + n_2) \end{array} \right.
\end{equation}
By solving 
for either $u$ or $v$, both of which can now be expressed linearly in $x$, and substituting back into the appropriate line of \eqref{double_eq_expl}, Diophantus obtains an equation for $x$ that is at worst quadratic. Moreover, the cases of \eqref{double_eq_expl} considered in the {\it Arithmetica} are such that, for a suitable choice of the parameter $\lambda$, Diophantus can ensure that the resulting equation for $x$ is linear (we will come back to the problem of the choice of $\lambda$ below). He thus finds a rational value for $x$, thereby solving the double equation. 

\subsection{Example} As an example, we give the double equation treated in the solution to Problem III.13, permitting ourselves the use of more than one unknown variable (the $x$ corresponds to Diophantus' unknown). The problem leads to the double equation
\begin{equation}
\label{III13_double_eq}
\left\{ ~~ \begin{array}{lcl} 4 x^2+15 x & = & u^2 \\ 4x^2-x-4 & = & v^2 \end{array} \right.
\end{equation}
With an eye towards the fact that
$$
(u+v)(u-v) = u^2 - v^2 = 16x+4 = (4x+1) \cdot 4,
$$
Diophantus sets
$$
u+v=4x+1,~u-v=4 ~~\Longrightarrow~~ u = 2x+\frac{5}{2},~ v = 2x-\frac{3}{2}.
$$ 
Substituting $u=2x+\dfrac{5}{2}$ into the top equation of \eqref{III13_double_eq} gives
\begin{equation}
\label{III17conclusion}
15x = 10x + \frac{25}{4} ~~\Longrightarrow ~~ x=\frac{5}{4}, ~u=5,~ v=1,
\end{equation}
and the equation is solved.

In the example above, it was clear that the factorization of $16x+4$ that we chose, which was $(4x+1)\cdot 4$, had to be such that the resulting equation \eqref{III17conclusion} became linear. We will discuss the appropriate choice of factorization in Section \ref{case_dist}. 

\subsection{Lines on quadrics} 
In this subsection and the next, we will summarize Diophantus' procedure, adding some geometric perspective. First, he subtracted the two equations making up his double equation \eqref{double_eq_expl}, as follows
\begin{equation}
\label{quadric}
u^2 - v^2 = (a_1-a_2) x^2 + (b_1 - b_2) x + (c_1 - c_2).
\end{equation}
Let $Q$ be the quadric surface in $\mathbb{P}^3$ defined by this equation. Diophantus subsequently relies on the fact that, in the cases which he considers, the right-hand side of \eqref{quadric} splits into linear factors, giving
\begin{equation}
\label{quadric_w_line}
(u+v)(u-v)=u^2 - v^2 = (m_1 x + n_1)(m_2 x + n_2),
\end{equation}
and uses this factorization to impose the following relations
\begin{equation}
\label{subst2}
\left\{ ~~ \begin{array}{lcl} u+v & = & \lambda \cdot (m_1 x + n_1) \\ u-v & = & \lambda^{-1} \cdot (m_2 x + n_2) \end{array} \right.
\end{equation}
for a suitable value of $\lambda$. 

Again, put in geometric terms, we can view \eqref{subst2} as defining a pencil of lines $\{L_\lambda\}_\lambda$ on the quadric $Q$ given by \eqref{quadric}; moreover, for rational values of $\lambda$, the line $L_\lambda$ is defined over the rationals. (It is well-known that smooth quadrics contain exactly {\it two} pencils of lines, or {\it rulings}; by interchanging the factors $m_1x+n_1$ and $m_2x+n_2$ in \eqref{subst2}, we obtain the other pencil.) 

In fact, the existence of the factorization \eqref{quadric_w_line} is {\it equivalent} to $Q$ containing a line defined over the rationals. Namely, suppose the right-hand side of \eqref{quadric} were irreducible, say with $\mathbb{Q}(\sqrt{d})$ as its splitting field. In that case, any automorphism $\sigma \in \operatorname{Aut}\overline{\mathbb{Q}}$ that interchanges $\pm \sqrt{d}$ also interchanges the two rulings of $Q$, as can be seen from \eqref{subst2}. Therefore $Q$ does not contain a line defined over the rationals. 

Let us note however that, conversely, the existence of the factorization \eqref{quadric_w_line} is not enough to construct a rational point on $C$. Indeed, take the system
\begin{equation}
\label{double_eq_fac_no_pt}
\left\{
\begin{array}{rcl}
3x^2 - 1 & = & u^2 \\
x^2 + 1 & = & v^2
\end{array}
\right.
\end{equation}
Here we have
$$
u^2 - v^2 = 2(x+1)(x-1),
$$
which defines a quadric with a line over the rationals. However, the smooth genus one curve $C$ defined by \eqref{double_eq_fac_no_pt} does not have rational points: there are no rational points at infinity, and the first equation has no $2$-adic (or $3$-adic) solutions.

Let us put translate this last fact back into more algebraic terms, so as to understand better the limitations of Diophantus' method: {\it having arrived at relations of the form \eqref{subst2}, it may not always be possible to choose $\lambda$ in such a way that a linear equation for $x$ is obtained.} We will now turn to the question of how Diophantus chooses his $\lambda$.

\subsection{The choice of the parameter $\lambda$}
\label{case_dist}

Solving \eqref{subst2}, we get for $u$ and $v$ the following linear expressions in $x$
$$
u = \frac{\lambda m_1 + \lambda^{-1} m_2}{2} x + \frac{\lambda n_1 + \lambda^{-1} n_2}{2}~, ~~ v = \frac{\lambda m_1 - \lambda^{-1} m_2}{2} x + \frac{\lambda n_1 - \lambda^{-1} n_2}{2}.
$$
For Diophantus' strategy to succeed, he needs to find a value for $\lambda$ such that substituting {\it either} the expression for $u$ into the first line of \eqref{double_eq_expl}, {\it or} substituting the expression for $v$ into the second line of \eqref{double_eq_expl}, causes either the constant or quadratic $x$-term to drop out. Now, for the double equations discussed in the {\it Arithmetica}, one of the following three conditions\footnote{As mentioned before, this subdivision is taken from Heath \cite[pp. 81--87]{heath}} is always satisfied:
\begin{itemize}
\item[(i)] $a_1=a_2$ is a (non-zero) square;
\item[(ii)] $a_1$ is square, and $a_2=0$;
\item[(iii)] $c_1=c_2=0$.
\end{itemize}
In every one of these cases, it is easy to show that a suitable value for $\lambda$ indeed exists. Among them, only case (iii) always leads to a curve of genus zero, and Diophantus' methods for dealing with it differ accordingly. We will therefore restrict our attention to cases (i) and (ii). 

\subsubsection{Case (i)} 
\label{cond_i}
First, let us assume condition (i) is satisfied. Put $a_1 = \alpha^2$. Since $u^2-v^2$ coincides with a linear polynomial in $x$ on $C$, one gets $m_1m_2=0$. If furthermore, as always in the {\it Arithmetica}, $m_1$ and $m_2$ are not {\it both} zero, one easily sees from the expressions for $u$ and $v$ given above that one can find $\lambda$ with the required property. Indeed, let us assume that $m_2=0$. Also, without loss of generality, we may suppose $n_2=1$, so that \eqref{subst2} becomes
\begin{equation}
\label{subst_cond_i}
\left\{ ~~ \begin{array}{lcl} u+v & = & \lambda \cdot (m_1 x + n_1) \\
 & = & \lambda \cdot \left((b_1-b_2)x + (c_1-c_2)\right) \\ u-v & = & \lambda^{-1}   \end{array} \right.
\end{equation}
Now consider the resulting expression for $u$
$$
u = \frac{\lambda m_1}{2}\, x + \frac{\lambda n_1 + \lambda^{-1} n_2}{2}.
$$
It is clear that $u^2$ is a quadratic polynomial in $x$ whose leading coefficient is a square, and that one can in fact choose $\lambda$ in such a way that the leading coefficient will equal $a_1$. Explicitly, one should choose
$$
\lambda = \dfrac{2\alpha}{m_1} = \dfrac{2\alpha}{b_1-b_2}.
$$

To put it in geometric terms: for the above value of $\lambda$, the line given by \eqref{subst_cond_i} intersects the point $(1:\alpha:\alpha:0)$ on $C$. It therefore intersects $C$ in a unique second rational point.

\subsubsection{Case (ii)}
\label{cond_ii}
Let us now start from condition (ii). We again put $a_1 = \alpha^2$. We scale the second equation to get $c_1=c_2$. Since $u^2-v^2$, when viewed as a polynomial in $x$, has no constant term, we find $n_1n_2$. As before, the case $n_1=n_2=0$ does not occur. Let us assume $n_2=0$. Without loss of generality, we may also set $m_1=m_2=\alpha$. It is again clear that we can get find a suitable $\lambda$; indeed, the system \eqref{subst2} becomes
\begin{equation}
\label{subst_cond_ii}
\left\{ ~~ \begin{array}{lcl} u+v & = & \lambda \cdot (m_1 x + n_1) \\
      & = & \lambda \cdot \left(\alpha x + (b_1-b_2)/\alpha \right) \\
 u-v & = & \lambda^{-1} \cdot \alpha x  \end{array} \right.
\end{equation}
We choose $\lambda = 1$. Now $v$ turns out to be a constant polynomial in $x$:
$$
v = \frac{\lambda m_1 - \lambda^{-1} m_2}{2} x + \frac{\lambda n_1}{2} = \frac{\lambda n_1}{2} = \dfrac{b_1-b_2}{2 \alpha}.
$$
By substituting this into the second line of \eqref{double_eq_expl}, we obtain an equation for $x$ that has no quadratic term.\footnote{Here, as in case (i), we do not claim that the resulting equation for $x$ is strictly of degree $1$: it might be that the linear term drops out as well. In that case, the solution must be interpreted as being ``at infinity''. Naturally, such cases do not occur in the {\it Arithmetica}.} 

Geometrically, we have intersected $C$ with the unique line from the pencil \eqref{subst_cond_ii} that passes through the point $(1:\alpha:0:0)$.

\subsubsection{Summary of the method}
The arguments from \ref{cond_i}--\ref{cond_ii} can be summarized as follows (thereby inevitably generalizing them): if $C$ contains a known rational point $P_1$, one chooses a rational value for $\lambda$ such that \eqref{subst2} defines a line $L$ containing $P_1$. Then $L$ will intersect $C$ in a unique second point $P$, which is rational.\footnote{It must be noted that the method as we have discussed it does not guarantee that the coordinates of $P$ are all positive. Indeed, this need not be the case. An example of such a failure is Problem V.2, where the double equation 
$$
x^2+20=u^2,~~4x+20=v^2
$$
leads to the equation $4x+20=4$, ``which is absurd'' (\textgreek{ὅπερ ἄτοπον}). (Likewise, one also cannot exclude the possibility that $P$ is a point at infinity, but this is more of an exceptional case, and it seems less likely that it is attested in the {\it Arithmetica}, although we have not made an exhaustive search.)} In algebraic terms, this means that the resulting equation is necessarily linear, or splits into linear factors.

In the rest of this section, we assume that $C$ is a (smooth) curve of genus one. We will explain the construction in terms of the group law on $C$, and describe an improvement suggested by Fermat, who realized that Diophantus' strategy can be applied iteratively.

\subsection{A contribution by Fermat}\label{fermat}
In his marginal commentary on Bachet's {\it Arithmetica}, Fermat wrote (translation by Heath \cite[p. 321]{heath})
\begin{quotation}
This double-equation gives, it is true, only one solution, but from this solution we can deduce another, from the second a third, and so on. In fact, when we have obtained one value for $x$ [say $x = a$], we substitute for $x$ in the equations the binomial expression consisting of $x$ {\it plus} the value found [i.e. $x+a$]. In this way we can find any number of successive solutions each derived from the preceding one.
\end{quotation}
We spend the rest of this section making Fermat's idea more explicit. It turns out that, in terms of the group law on $C$ (considered as an elliptic curve), Fermat roughly proposes a repeated {\it triplication} process, starting from an initial solution $P$. This means that the numbers involved in writing down the subsequent solutions would grow {\it very} fast.

We return to the general setup of double equations, that is to the system
\begin{equation}
\label{double_eq_expl2}
\left\{ ~~ \begin{array}{lcl} a_1 x^2+b_1 x+c_1 & = & u^2 \\ a_2x^2+b_2x+c_2 & = & v^2 \end{array} \right.
\end{equation}
As before, let us write $C$ for the curve defined by this system, and we assume that $C$ is smooth and irreducible, and hence of genus one. Let us also assume that $C$ has a rational point. Without loss of generality, we may assume that $C$ has a rational point at infinity, which is equivalent to saying that
$$
a_1 = \alpha_1^2,~~ a_2 = \alpha_2^2
$$
with $\alpha_1$ and $\alpha_2$ rational numbers. We write
\begin{align*}
P_1 & = (1:\alpha_1:\alpha_2:0),\\\
P_2 & = (1:-\alpha_1:-\alpha_2:0),\\
P_3 & = (1:-\alpha_1:\alpha_2:0),\\
P_4 & = (1:\alpha_1:-\alpha_2:0)
\end{align*}
for the points at infinity\footnote{One could wonder about the relations between the elements $(P_i)-(P_j)$ of the Jacobian $\operatorname{Jac}(C)$ of $C$, which parametrizes linear equivalence classes of divisors of degree zero on $C$. In general, it turns out that $(P_2)-(P_1)$ is of order $2$, that $(P_4) - (P_3) = (P_2) - (P_1)$, and that $(P_3)-(P_1)$ is generically of infinite order. For the last statement, see Section \ref{double_eqs_infinite_order}.}. Since $C$ is smooth, $\alpha_1$ and $\alpha_2$ are not both zero, so there are either $2$ or $4$ points at infinity.

In order to obtain a quadric $Q$ with a rational line and containing $C$, we write\footnote{This step is to be compared to Section \ref{case_dist}, case (i). Obviously, taking the linear combination \eqref{gen_lin_comb} of the two equations of \eqref{double_eq_expl2} is equivalent to first scaling both equations by $a_2$ and $a_1$ respectively and proceeding as in Section \ref{case_dist}.}
\begin{equation}
\label{gen_lin_comb}
a_2 u^2 - a_1 v^2 = (\alpha_2 u + \alpha_1 v)(\alpha_2 u - \alpha_1 v) = m x + n, 
\end{equation}
with
$$
m = a_2 b_1 - a_1 b_2, ~ n = a_2 c_1 - a_1 c_2,
$$
which leads us to consider the line $L_\lambda$ on $Q$ given by
\begin{equation*}
\label{subst3}
\left\{ ~~ \begin{array}{lcl} \alpha_2 u + \alpha_1 v & = & \lambda \cdot (m x + n) \\ \alpha_2 u - \alpha_1 v & = & \lambda^{-1} \end{array} \right.
\end{equation*}

Note that the line $L_0$ passes through $P_3$ and $P_4$. Following Diophantus' strategy, we would take the unique line $L_{\lambda}$ passing through $P_1$, and define $P$ as the second point of intersection between $L_{\lambda}$ and $C$.

Given the newly found point $P$, it is Fermat's proposal to replace $x$ by $x'+x(R)$, so that, he seems to imply, the original procedure may be applied afresh. Indeed, the effect of Fermat's substitution is that we end up with a new system, say defining a curve $C'$.
\begin{equation*}
\left\{ ~~ \begin{array}{lcl} a'_1 x'^2+b'_1 x'+c'_1 & = & u^2 \\ a'_2 x'^2+b'_2x'+c'_2 & = & v^2 \end{array} \right.
\end{equation*}
Since $x'=0$ corresponds to a rational point on $C'$, we now have that $c_1'$ and $c_2'$ are both squares, and it is quite clear that we may proceed analogously to before, except that the role of $x$ is now taken by $1/x'$ rather than $x'$ itself.  

However, rather than having to change the curve $C$ each time a new point is found, it is of course more convenient to express Fermat's idea in terms of the original curve. Let $P=(x_0,u_0,v_0)$ be any point on $C$. We will further need to consider the points
$$
P' = (x_0,-u_0,v_0),~~P''=(x_0,u_0,-v_0).
$$
Completely analogously to before, we can construct a quadric $Q_P$ associated to $P$, and which contains a line passing through $P'$ and $P''$. We simply define $Q_P$ by the equation
\begin{equation}
\label{more_gen_lin_comb}
v_0^2 u^2 - u_0^2 v^2 = (v_0 u + u_0 v)(v_0 u - u_0 v) = \ell x^2 + m x + n, 
\end{equation}
with
$$
\ell = v_0^2 a_1 - u_0^2 a_2, ~ m = v_0^2 b_1 - u_0^2 b_2, ~ n= v_0^2 b_1 - u_0^2 b_2.
$$
Since $(x_0,u_0,v_0)$ lies on $C$, and therefore on $Q_P$, the right-hand side of \eqref{more_gen_lin_comb} is reducible, having $x-x_0$ as a factor. Write the right-hand side as $(x-x_0) h(x)$ with $h$ of degree $\leq 1$, so that $Q_P$ is given by
$$
(v_0 u + u_0 v)(v_0 u - u_0 v) = (x-x_0) h(x).
$$
Then on $Q_P$ we have a pencil of lines $L_{P,\lambda}$, given by
\begin{equation}
\label{subst4}
\left\{ ~~ \begin{array}{lcl} v_0 u + u_0 v & = & \lambda \cdot h(x) \\ v_0 u - u_0 v & = & \lambda^{-1} \cdot (x - x_0) \end{array} \right.
\end{equation}
As can be seen from the equations, the line $L_{P,0}$ again passes through $P'$ and $P''$. If we denote by $L_{P,\lambda}$ the line through $P$, then  the second point of intersection $R$ between $L_{P,\lambda}$ and $C$ would again be the point we would find by applying the idea of Diophantus--Fermat.

\subsubsection{The group law on $C$}
We continue with the notation established so far. Now, since $L_{P,0}$ intersects $C$ in $P'$ and $P''$, and $L_{P,\lambda}$ intersects $C$ in $P$ and $R$, and since $L_{P,0}$ and $L_{P,\lambda}$ belong to the same linear system, we have the following linear equivalence of divisors on $C$
$$
P + R \sim P' + P''.
$$

We now want to express $R$ solely in terms of $P$ and points on $C$ that were chosen independently of $P$. To do this, we claim that we have the following linear equivalences between divisors on $C$
$$
P' \sim P_1 + P_3 - P, ~~ P'' \sim P_1 + P_4 - P. 
$$
The proof of the claim involves some theory about elliptic curves (as can be found in \cite{silverman}). We prove the first claim; the proof of the second is analogous. Let $\phi$ be the automorphism of $C$ sending a point $(x,u,v)$ to $(x,-u,v)$. If we consider $C$ as an elliptic curve with identity $P_1$, then it follows from the theory of elliptic curves that $\phi$ can be obtained by composing an endomorphism $\epsilon$ of $C$ with a translation by some point $T_0$, that is, we have
$$
\phi(T) = \epsilon T + T_0, 
$$
for all points $T$ on $C$. Since $\phi(P_1) = P_3$, we know that $T_0 = P_3$. By the fact that $\phi \circ \phi$ is the identity on $C$, the endomorphism $\epsilon$ must be an automorphism. Since the generic $C$ of the form \eqref{double_eq_expl2} does not have any automorphisms (over the algebraic closure of $\mathbb{Q}$) other than $\pm 1$, we have $\epsilon = \pm 1$. Suppose first that $\epsilon = 1$. Again arguing from the generic case, this is impossible, since $P_3$ is generically of infinite order, which contradicts $\phi \circ \phi = \operatorname{id}$. Therefore $\epsilon = -1$, and 
$$
\phi(T) = P_3 - T, 
$$
for all points $T$ on $C$. Putting this in terms of linear equivalences of divisors, we obtain the claim.\footnote{Incidentally, the claim just proven also implies that $P_2$ is of order $2$ on $C$ and that $P_2+P_3=P_4$ (still taking $P_1$ as the identity). Let $\psi$ be the automorphism of $C$ sending $(x,u,v)$ to $(x,u,-v)$, and $\chi$ the automorphism sending $(x,u,v)$ to $(x,-u,-v)$, so that $\chi = \psi \circ \phi$, with $\phi$ as above. Applying the claim, we find that $\psi(T)=P_4-T$ for all points $T$ on $C$, so that  $\chi = \psi \circ \phi$ is translation by $P_4-P_3$, which must be of order $2$ since $\chi \circ \chi$ is the identity on $C$. Moreover, since $\chi(P_1)=P_2$, we find that $P_4-P_3=P_2$, which proves what we wanted.}

Having established the claim, we can finally express $R$ in terms of $P$, as follows
$$
R \sim 2P_1 + P_3 + P_4 - 3 P 
$$
or, in terms of the group law with $P_1$ as the identity
$$
R = P_3 + P_4 - 3P.
$$
If we take $P=P_1$ as our zero-th solution, then $R = P_3+P_4$ will be the first; the next few iterations are $-2R$, $7R$, $-20R$, $61R$, $\ldots$, and  in general, the $n$-th solution will be
$$
\dfrac{1-(-3)^n}{4} \cdot R.
$$

\section{Problem III.17: a K3 surface of degree $8$}
\label{secIII17}

The seventeenth problem of Book III runs as follows:
\begin{quotation}
To find two numbers such that their product added to both or to either gives a square.
\end{quotation}

\subsection{Diophantus' solution of Problem III.17}

Let the two numbers be $x$ and $y$. Diophantus chooses $y=4x-1$, so that the product $x(4x-1)$ is $x$ less than the square $4x^2$, which means that one condition is satisfied. The two remaining conditions give the double equation (here explained in section \ref{double_eqs})
$$
\left\{ \begin{array}{lcl}
 4x^2 + 4x - 1 & = & u^2 \\
 4x^2 + 3x - 1 & = & w^2
 \end{array} \right.
$$
In particular he obtains
$$
(u+w)(u-w) = u^2 - w^2 = x.
$$
Following the methods for treating double equations (see Section \ref{double_eqs}), which he assumes the reader knows by this point, Diophantus sets $u + w = 4x$ and $u - w = \dfrac{1}{4}$, so that $w = 2x - \dfrac{1}{8}$. Substituting in the second equation yields
$$
 4x^2 + 3x - 1 = \left(2x - \dfrac{1}{8}\right)^2 = 4x^2 - \dfrac{x}{2} + \dfrac{1}{64} 
$$
which gives $\dfrac{7x}{2}=\dfrac{65}{64}$, and so $x = \dfrac{65}{224}$, and $y = 4x-1 = \dfrac{9}{56}$.

\subsection{The geometry of Problem III.17}

The problem leads to a surface $X$ given by 
\begin{equation}
\label{III17a}
\left\{ ~~ \begin{array}{lcl} xy +x+ y & = & u^2 \\ xy + x & = & v^2 \\ xy +y & = & w^2 \end{array} \right.
\end{equation}

These equations define a complete intersection of three quadrics in $\mathbb{P}^5$ that is singular, but with at worst rational double points as singularities (calculated with \texttt{magma} \cite{magma}), which means that $X$ is birationally equivalent\footnote{\label{K3def}The usual definition of a \defi{K3 surface} is a ``projective, smooth, and geometrically connected surface $Y$ such that the canonical sheaf $\omega_Y$ is trivial and $\operatorname{H}^1(Y,\mathscr{O}_Y)=0$''. However, here and in the next section, we will mostly employ the term K3 surface in a broader sense, and use it to refer to a surface that properly speaking is only birationally equivalent to one. Since the study of rational points on projective non-singular varieties is largely insensitive to switching between birational models, and the methods used by both Diophantus and us are mostly birational in nature, there seems to be some justification for this from our point of view.} to a K3 surface (see Theorem \ref{singularK3} in Section \ref{singK3}). Diophantus' solution leads to a pencil of curves whose generic member is a smooth intersection of two quadrics in $\mathbb{P}^3$, hence a curve of genus $1$. On this pencil Diophantus constructs a rational section as follows. 

Diophantus chooses $y=4x-1$, which we can naturally generalize to $y = t^2x - 1$, with $t=2$ being Diophantus' choice. Then the second equation reads
$$
xy + x  = x(t^2x-1)+x = t^2x^2 =  v^2,
$$
which is fulfilled for $v=\pm t x$, so that we may drop the second equation from our considerations. Reversing this logic, if we define $t=v/x$, then the second equation implies 
$$
x(y+1)= v^2 ~~ \Longrightarrow ~~ y=\dfrac{v^2}{x}-1 = t^2x - 1.
$$
Therefore, in terms of the rational map
\begin{align*}
\pi \colon X & \dashrightarrow \mathbb{P}^1 \\\
(x,y,u,v,w) & \longmapsto v/x
\end{align*}
the argument put forward by Diophantus comes down to restricting attention to the fibre $X_{t_0}$ of $\pi$ over a rational point $(t_0:1)$ on $\mathbb{P}^1$.

We are then left with the following set of equations in $u,w,$ and $x$:
\begin{equation}
\label{III17fibre}
X_t : ~ \left\{ ~~ \begin{array}{lcl} t^2 x^2 + t^2 x - 1 & = & u^2 \\ t^2 x^2 + (t^2 - 1)x - 1 & = & w^2 \end{array} \right.
\end{equation}
or written projectively, 
\begin{equation}
\label{III17fibre_hom}
X_t : ~ \left\{ ~~ \begin{array}{lcl} t^2 Y^2 + t^2 YZ - Z^2 & = & U^2 \\ t^2 Y^2 + (t^2 - 1)YZ - Z^2 & = & W^2 \end{array} \right.
\end{equation}
It is easy to see that, for generic $t$, these equations define a smooth intersection of two quadrics in $\mathbb{P}^3$. This follows for example from the fact that the polynomial
$$
(t^2 x^2 + t^2 x - 1)(t^2 x^2 + (t^2 - 1)x - 1)
$$
has four distinct zeros for almost all $t$ (see footnote \ref{smoothness_crit}). Hence the generic fibre of $\pi$ is a smooth curve of genus $1$. Since $\pi$ has natural sections $(Y:U:W:Z)=(1:\pm t: \pm t:0)$, we may select one of them as the origin of the group law on $X_t$. We choose to let $\mathcal{O}=(1:t:t:0)$ be the identity on $X_t$.

Subtracting the second equation from the first, we get
$$
x = u^2 - w^2 = (u-w)(u+w).
$$
By Diophantus' standard methods for dealing with ``double equations'' (see Section \ref{double_eqs}), we put $u+w=2tx$ and $u-w=1/(2t)$, we get 
$$
u = tx + \frac{1}{4t},~~ w = tx - \frac{1}{4t}.
$$
Substituting the expression for $w$ in the bottom equation of \eqref{III17fibre} gives
$$
(t^2 - 1)x - 1 = - \frac{1}{2} x  +  \frac{1}{16t^2}  ~~ \Longrightarrow ~~ x = \frac{16t^2 + 1}{16t^4-8t^2},
$$
which allows to solve for $y,u,v,w$, and gives the map $\sigma\colon \mathbb{P}^1 \rightarrow X$ with $x$ as above and $y,u,v,w$ defined  by
$$
y=\frac{9}{16t^2-8}, ~~ u=\frac{20t^2 - 1}{16t^3-8t},~~v=\frac{16t^2 + 1}{16t^3-8t},~~w=\frac{12t^2 + 3}{16t^3-8t}.
$$
In the original solution, Diophantus chose $t=2$, which leads to the point
$$
P=(x,y,u,v,w)=\left( \dfrac{65}{224},\dfrac{9}{56},\dfrac{79}{112},\dfrac{65}{112},\dfrac{51}{112} \right).
$$

\subsection{The order of $P$ on $X_2$} 
\label{double_eqs_infinite_order}

It is natural to ask whether the {\it order} of $P$ on the fibre $X_2$ is infinite. We shall prove now that it is. (This subsection is somewhat more technical than most of this article, since it uses some facts about the arithmetic of elliptic curves. The necessary background can be found in \cite{silverman} and \cite[Chapter IV]{silverman2}.)

First, a straightforward application of the smoothness criterion mentioned in footnote \ref{smoothness_crit} gives that $X_2$ has good reduction modulo $7$, in the strong sense that \eqref{III17fibre_hom} defines a smooth genus $1$ curve $\widetilde{X}_2$ over $\mathbb{F}_7$. Denoting reduction modulo $7$ of both integers and points on $X_2$ by vertical bars, we have (using projective coordinates)
$$
\overline{P} = (\overline{65}:\overline{158}:\overline{102}:\overline{224}) = (\overline{2}:\overline{4}:\overline{4}:\overline{0}) = \overline{\mathcal{O}_2},
$$
so $P$ is contained in the kernel $X_2^1(\mathbb{Q}_7)$ of the reduction homomorphism
$$
X_2(\mathbb{Q}_7) \rightarrow \widetilde{X}_2(\mathbb{F}_7).
$$
Since $X_2$ has good reduction modulo $7$, the subgroup $X_2^1(\mathbb{Q}_7)$ is a formal group over $\mathbb{Z}_7$, and therefore torsion-free (see \cite[IV.6 and VII.3]{silverman}). 

This argument also establishes that the point $P$ constructed in Section \ref{case_dist} is generically of infinite order. 

\subsection{Summary} 
Starting out with the K3 surface $X$ defined by \eqref{III17a}, we have constructed a rational map
\begin{align*}
\pi \colon X & \dashrightarrow \mathbb{P}^1 \\\
(x,y,u,v,w) & \longmapsto v/x
\end{align*}
which seemed to offer most natural way of viewing Diophantus' solution. The generic fibre of $\pi$ is a smooth curve of genus $1$. As we have seen following Diophantus' method, the map $\pi$ has a section $\sigma$. In fact, the surface $X$ is an \defi{elliptic K3 surface}: instead of $\sigma$, there are also the obvious sections at infinity (these would have had little meaning for Diophantus, of course); we chose $\mathcal{O}$ to be one of these. The section $\sigma$, together with the calculation from the previous subsection, proves that the rank of the elliptic surface $(X,\mathcal{O})$ must be at least $1$.

\section{Problem IV.18: a K3 surface of degree $6$}
\label{secIV18}

Here is the statement of the eighteenth problem of Book IV:
\begin{quotation}
To find two numbers such that the cube of the first added to the second gives a cube, and the square of the second added to the first gives a square.
\end{quotation}

\subsection{Diophantus' solution to Problem IV.18}

Letting $x$ be the first number and $y$ the second, Diophantus proposes to make the sum $x^3+y$ equal to $8$, so that the first requirement is met, and he has $y=8-x^3$. Substituting this into the second equation, he obtains
$$
x + (8-x^3)^2 = v^2 ~~ \Longrightarrow  ~~ v^2 = x^6 - 16x^3 + x + 64.
$$
He now chooses $v = x^3+8$, and he obtains
$$
(x^3+8)^2 = x^6+16x^3 + 64 = x^6 - 16x^3 + x + 64 ~~\Longrightarrow~~ 32x^3 = x,
$$
which has no positive rational solutions. Hence he observes that the coefficient $32$ arose as $16+16$, and both occurrences of $16$ arose as $2 \times 8$, which is the same $8$ as was chosen in the beginning. Hence, Diophantus says, instead of the cube $8$ we should have chosen a cube $z^3$ such that $4z^3$ is a square, say $16z^2$, so that $z=4$. Replacing his original choice $8$ by $z^3=64$ then, Diophantus obtains 
$$
x + (64-x^3)^2 = v^2 ~~ \Longrightarrow  ~~ v^2 = x^6 - 128x^3 + x + 4096.
$$
Choosing $v=x^3+64$, he then ends up with
$$
(x^3+64)^2 = x^6+128x^3 + 4096 = x^6 - 128x^3 + x + 4096 ~~\Longrightarrow ~~256x^3 = x,
$$
hence $x = \dfrac{1}{16}$ and $y = 64 - x^3 = \dfrac{262143}{4096}$.
\subsection{The geometry behind Diophantus' solution}

The equations are as follows:
\begin{equation}
\label{IV18a}
\left\{ ~~ \begin{array}{lcl} x^3 + y & = & u^3 \\ x + y^2 & = & v^2 \end{array} \right.
\end{equation}

These define a singular complete intersection $X$ of a quadric and a cubic hypersurface in $\mathbb{P}^4$, whose singularities are all rational double points (calculated with \texttt{magma} \cite{magma}). Hence $X$ is a K3 surface (see Section \ref{singK3} and footnote \ref{K3def} for our use of the term ``K3 surface''). Following Diophantus' solution, we will end up with a rational curve on $X$.

Diophantus imposes the equation $x^3+y=u_0^3$ for a fixed $u_0$, which translated into modern terms means examining the fibres of the map
\begin{align*}
\pi \colon X & \dashrightarrow \mathbb{P}^1 \\\
(x,y,u,v) & \longmapsto u
\end{align*}
Thus one writes $y=u_0^3-x^3$, eliminating $y$, and use the second equation to get
$$
x + (u_0^3-x^3)^2 = v^2 ~~ \Longrightarrow  ~~ v^2 = x^6 - 2 u_0^3 x^3 + x + u_0^6.
$$
This gives a genus-$2$ hyperelliptic curve for generic values of $u_0$ (for example for $u_0=0$, as is easy to verify), but without an obvious rational point in general. Diophantus does attempt the substitution $v=x^3+u_0^3$, which eventually yields
$$
4 u_0^3 x^3 = x,
$$
which fails, because for general $u_0$ the coefficient in the left-hand side is not a square.

Of course, the next step is to remedy the failure by choosing $u_0$ so that $4u_0^3$ is square. To us, it would be natural to consider this condition in terms of the prime factorization of $u_0$, which immediately yields $u_0$ being a square as an equivalent condition. However, judging from the solution above, this approach is less natural for Diophantus, who apparently prefers to think algebraically rather than arithmetically.\footnote{In the four ``Arabic'' books of the {\it Arithmetica} in particular, one encounters a fair number of problems whose solutions could have been shortened by an appeal to unique prime factorization (of which Euclid was already aware; see Elements VII, Propositions 30--32).} 

Instead, Diophantus assumes for $4u_0^3$ the value $16u_0^2$, which for $u_0$ yields the linear equation $4u_0=16$, so that $u_0=4$. By varying the coefficient $16$, Diophantus could of course have obtained any square value of $u_0$. The solution to this little subproblem can be reinterpreted as parametrizing the rational points on the cuspidal curve
$$
4u^3=w^2,
$$
a theme that recurs a number of times in Diophantus' work, especially in the (Arabic) Book 4.

The idea to restrict to square values of $u_0$ leads us to consider the base-change of the fibration $\pi$ by the map 
\begin{align*}
\phi \colon \mathbb{P}^1 & \rightarrow \mathbb{P}^1 \\\
t & \mapsto t^2
\end{align*}
from the $t$-line to the $u$-line. 

We denote by $X'$ the base-change of $X$ by $\phi$. It is not hard to see that $X'$ is geometrically irreducible. For example, consider the curve $C \subset X$ defined by $x=0$, which is defined by
$$
v^2 = u^6+u.
$$
The function field of $C$ is $K=\mathbb{Q}(u)[\sqrt{u^6+u}]$, and it is an elementary exercise in field theory that $\sqrt{u} \notin K$, and therefore $u$ is  not a square in the function field of $X$. 

The equations of $X'$ are now given by
\begin{equation}
\label{IV18b}
\left\{ ~~ \begin{array}{lcl} x^3 + y & = & t^6 \\ x + y^2 & = & v^2 \end{array} \right.
\end{equation}
that is, we have simply replaced $u$ in \eqref{IV18a} by $t^2$. Proceeding as before, we eliminate $y$ and substitute $y=t^6-x^3$ into the second equation, and we obtain
$$
x + (t^6-x^3)^2 = v^2 ~~ \Longrightarrow  ~~ v^2 = x^6 - 2t^6x^3 + x + t^{12}.
$$
Analogously to before, we choose $v = x^3+t^6$, and obtain
$$
(x^3+t^6)^2 = x^6+2t^6 x^3 + t^{12} = x^6 - 2t^6 x^3 + x + t^{12} ~~\Longrightarrow~~ 4t^6x^3 = x,
$$
which gives $x=\dfrac{1}{2t^3}$, and so $y=t^6-x^3=\dfrac{8t^{15}-1}{8t^9}$ (and $v=\dfrac{8t^{15}+1}{8t^9}$). 

\subsection{Summary} 
The K3 surface $X$ given by \eqref{IV18a} has a fibration
$$
\pi \colon X \dashrightarrow \mathbb{P}^1
$$
into curves of genus $2$, but without an obvious section. After base-changing $\pi$ along the map
\begin{align*}
\phi \colon \mathbb{P}^1 & \rightarrow \mathbb{P}^1 \\\
t & \mapsto t^2
\end{align*}
we obtain a new surface $X'$, defined by \eqref{IV18b}, together with a fibration $\pi' \colon X' \dashrightarrow \mathbb{P}^1$. Since $\pi'$ is obtained by base-change from $\pi$, the generic fibre is again a curve of genus $2$, and moreover $\pi'$ has a section, whose image on $X$ is birational to $\mathbb{P}^1$.

\section{Problem IV.32: a del Pezzo surface of degree $4$}

The thirty-second problem of Book IV runs thus:

\begin{quotation}
To divide a given number into three parts, such that the product of the first and second, after the third is either added or subtracted, yields a square.
\end{quotation}

\subsection{Diophantus' solution to Problem IV.32}

We continue with a translation of Diophantus' solution.

\begin{quotation}
Let the given number be $6$. Let the third part be $x$, and the second a fixed value less than $6$, say $2$, so that the first part is $4-x$. The problem is now to find a value for $x$ such that  the two numbers
$$
2 \cdot (4 - x) \pm x
$$
are both squares, which yields the double equation
\begin{equation}
\label{IV32provisional}
\left\{ ~~ \begin{array}{lcl} 8 - x & = & \textnormal{square} \\ 8-3x & = & \textnormal{square} \end{array} \right.
\end{equation}

We could have solved this system if the ratio between the linear coefficients had been a square.\footnote{In fact, Diophantus remarks of the problem that he ended up with, that it \textgreek{οὐ	ῥητόν ἐστι}, i.e. ``it is not rational'', which he says elsewhere of quadratic equations in one variable with no rational root (e.g.~in IV.31), or of the roots of such equations (e.g.~in IV.9). As it happens, the system under consideration does not have a rational solution either, but Diophantus does not show this.} 

[Therefore we should choose a value for the second part, instead of $2$, that does give a square ratio between the linear coefficients.]

Now the coefficient of $x$ in the first equation is one less than $2$ [which was our choice for the second part], and the coefficient of $x$ in the second equation is one more than $2$. So we are led to the problem of finding some number, instead of $2$, such that this number plus one has to the same number minus one the ratio of a square number to a square number.

Let the value of the second part be $x$.\footnote{Of course, Diophantus only has one symbol to denote an unknown quantity. We follow him in this with our double use of the symbol $x$.} Now the linear coefficients that will result from this choice are equal to $x+1$ and $x-1$. We want that their ratio is that of a square number to a square number, say that of $4$ to $1$. Now $4$ times $x-1$ equals $4x-4$, and $1$ times $x+1$ equals $x+1$. So $4x-4=x+1$, and $x$ becomes $\dfrac{5}{3}$. 

Therefore we now take the second part to be $\dfrac{5}{3}$, whereas we denote the third part by $x$, which means that the first part will be $\dfrac{13}{3} - x$. We therefore want
$$
\frac{65}{9}-\frac{5x}{3} \pm x 
$$
to  both be squares, which yields the system
\begin{equation*}
\left\{ ~~ \begin{array}{lcl} 65 - 6x & = & \textnormal{square} \\ 65-24x & = & \textnormal{square} \end{array} \right.
\end{equation*}
or equivalently
\begin{equation*}
\left\{ ~~ \begin{array}{lcl} 260 - 24x & = & \textnormal{square} \\ 65-24x & = & \textnormal{square} \end{array} \right.
\end{equation*}
Given a system such as this\footnote{This type of problem would be familiar to the reader of the {\it Arithmetica} at this point: in Problem II.11, Diophantus had already dealt with double equations such as the present system (see Section \ref{double_eqs} for the term ``double equation'' and for a full discussion of Problem II.11).
}, one proceeds as follows. We write the difference $260-65=195$ between the equations as a product of two factors, say $195 = 15 \times 13$, and take the first square to be half the sum of the factors squared\footnote{It is of course clear what happens here: the system $u^2-v^2=ab$ is rewritten as $(u+v)(u-v)=ab$, whereupon Diophantus chooses $u+v=a$ and $u-v=b$, which leaves him with linear equations for $u$ and $v$.}, and the second to be half the difference squared, that is, the squares will be $14^2=196$ and $1^2=1$. This finally yields $x=\dfrac{8}{3}$. Hence the first part is $\dfrac{5}{3}$, the second is $\dfrac{5}{3}$, and the third is $\dfrac{8}{3}$.
\end{quotation}

\subsection{The geometry behind Diophantus' solution to IV.32}
Diophantus asks, given a rational number $n$, for rational numbers $x,y,z$, such that there exist further rational numbers $u,v$ such that
\begin{equation}
\label{IV32a}
\left\{ ~~ \begin{array}{lcl} x + y + z & = & n \\ xy - z & = & u^2 \\ xy + z & = & v^2 \end{array} \right.
\end{equation}
Eliminating $z$ by writing $z=n-x-y$, we find that Problem IV.32 again leads to an intersection of two quadrics in $\mathbb{P}^4$:
\begin{equation}
\label{IV32b}
\left\{ ~~ \begin{array}{lcl} xy +x +y -n & = & u^2 \\ xy -x-y +n & = & v^2 \end{array} \right.
\end{equation}
which is  smooth for generic values of $n$, hence a del Pezzo surface $X$ of degree $4$ over $\mathbb{Q}$ by \cite[Proposition IV.16]{beauville}.

As is his custom, Diophantus starts out his solution by choosing a particular value for $y=y_0$ (he chooses $y_0=2$), which leaves him with the system
\begin{equation*}
\left\{ ~~ \begin{array}{lcl} y_0(n - y_0) + (1-y_0)z & = & u^2 \\ y_0(n - y_0) - (1+y_0)z & = & v^2 \end{array} \right.
\end{equation*}
Upon eliminating $z$, we obtain the conic
\begin{equation}
\label{IV32insolubleconic}
(y_0+1)u^2 - (y_0-1)v^2 + 2y_0(y_0-n) = 0.
\end{equation}
However, Diophantus finds that he cannot solve this equation for arbitrary $y_0$. Indeed, if we go along with him and take $n=6$ and $y_0=2$, then \eqref{IV32insolubleconic} becomes
\begin{equation}
\label{IV32insolubleconicspecialized}
3u^2 - v^2 - 16 = 0
\end{equation}
which has no solution in rational numbers, as can be established by considering solutions over the fields of $2$-adic and $3$-adic numbers.

The substitution $y=y_0$ leads to a conic bundle structure 
\begin{align*}
\pi \colon X & \dashrightarrow \mathbb{P}^1 \\\
(x,y,u,v) & \longmapsto y
\end{align*}
However, as the insolubility of \eqref{IV32insolubleconicspecialized} shows, this time $\pi$ does not have a section. However, Diophantus changes things around by assuming that $y_0$ is such that
\begin{equation}
\label{IV32pullback}
\frac{y_0+1}{y_0-1}=\textnormal{square}=t_0^2~~\textnormal{(say)}~\Longrightarrow ~  y_0 = \frac{t_0^2+1}{t_0^2-1}.
\end{equation}
In fact, he takes $t_0=2$, so that $y_0=\dfrac{5}{3}$. He proceeds by substituting $y_0=(t_0^2+1)/(t_0^2-1)$ into his equations, and after clearing denominators he ends up with the system
\begin{equation}
\label{IV32fibre}
\left\{ ~~ \begin{array}{lcl} (t_0^2+1)((n-1)t_0^2 - n-1) - 2(t_0^2-1) z & = & u^2 \\ (t_0^2+1)((n-1)t_0^2 - n-1) - 2t_0^2(t_0^2-1) z & = & v^2 \end{array} \right.
\end{equation}
(He of course has $t_0=2$ and $n=6$.) He then multiplies the first equation by $t_0^2$, subtracts the second equation from the first, and ends up with the conic 
\begin{equation}
\label{IV32genericconic}
t_0^2 u^2 - v^2 = (t_0^4-1)((n-1)t_0^2-n-1)
\end{equation}
To find solutions to this equation, Diophantus effectively proceeds as follows.\footnote{Of course, we could have inferred the existence of a parametrization of the rational solutions to \eqref{IV32genericconic} from the fact that the projective conic defined by \eqref{IV32genericconic} has a rational point at infinity.} First, he factors the right-hand side as $\lambda_0 \cdot \mu_0$ (he takes $\lambda_0 = 15, \mu_0 = 13$). Then since
$$
(t_0u + v)(t_0u - v) = t_0^2 u^2 - v^2 = \lambda_0 \cdot \mu_0
$$
and puts $t_0u + v = \lambda_0$, $t_0u - v = \mu_0$, which then leads to
\begin{equation}
\label{IV32paramfibre}
t_0u = \frac{\lambda_0 + \mu_0}{2}, ~ v = \frac{\lambda_0 - \mu_0}{2}, 
\end{equation}
and this finally allows us to solve for $z$ via \eqref{IV32fibre}.

\subsection{Summary}
As discussed, the map $\pi$ endows $X$ with the structure of a conic bundle; however, $\pi$ does not have a section. However, the substitution for $y_0$ as given by \eqref{IV32pullback} leads one to introduce the following rational map
\begin{align*}
\phi \colon \mathbb{P}^1 & \rightarrow \mathbb{P}^1 \\\
t & \mapsto \frac{t^2+1}{t^2-1}
\end{align*}
and to consider the base-change $\pi'$ of $\pi$ by the map $\phi$. In fact, the family of conics \eqref{IV32fibre}, where $t_0$ ranges over the rationals, form only a proper subset of the fibres of $\pi$ over rational points on $\mathbb{P}^1$, but they correspond exactly to the rational fibres of $\pi'$. Now the parametrization \eqref{IV32paramfibre}, which essentially followed from Diophantus' argument, shows that the conic fibration
$$
\pi' \colon X' \dashrightarrow \mathbb{P}^1,
$$
has a section.\footnote{We should in fact check that $X'$ is geometrically integral, i.e.~that the function $f=(y+1)/(y-1)$ ``is not a square already'' in the function field of $X$. To this end, consider the curve $C$ defined by $x=0$; its equations are $y-n=u^2$, $-y+n=v^2$, so that the function field of $C$ is $K=\mathbb{Q}(y)[\sqrt{-1},\sqrt{y-n}]$. It follows easily from basic field theory that $(y+1)/(y-1)$ is not a square in $K$.} Hence $X'$ is birational to $\mathbb{P}^2$. Moreover, $X'$ dominates $X$, so we find that $X$ is dominated by a surface that is birationally equivalent to $\mathbb{P}^2$. This gives the parametrization of $X$ from above.


\section{Problem V.29: a del Pezzo surface of degree $2$}
\label{secV29}

Problem V.29 runs as follows:
\begin{quotation}
To find three squares such that the sum of their squares is a square.
\end{quotation}

\subsection{Diophantus' solution to Problem V.29}
Diophantus proposes to take $y=2$ and $z=3$. He then gets $x^4+97 = w^2$. Taking the substitution $w=x^2-10$, so that $w^2=x^4-20x^2+100$, he ends up with $20x^2=3$, which does not give a rational solution. Hence, he says, we are looking for rational numbers $p,q,$ and $m$, such that 
$$
\frac{m^2-p^4-q^4}{2m}
$$
is a square.\footnote{Diophantus does not explain this further, so I will add the argument here. He wants to choose values $y=p$ and $z=q$ such that his trick for cancelling the $x^4$ term results in a polynomial equation with a rational solution; in order for this to be possible, there must exist $m$ such that expanding
$$
x^4-p^4-q^4=(x^2-m)^2
$$
gives a quadratic equation with rational roots. This yields the criterion above.}${}^{\textnormal{,}}$\footnote{It is interesting to give a more literal translation of this passage, so that we may see Diophantus at work at a moment where the algebra is particularly involved from his perspective. After obtaining $20x^2=3$, he goes on: ``Now if both [coefficients] had been squares, the unknown would have been found. So it is reduced to finding two squares and a certain number such that the square of this [last-mentioned number] minus the squares of the searched-for [squares] would make some number which has the same ratio  to twice the number from the beginning as that of a square number to a square number.'' That is, in our notation from above, the ratio between $m^2 - p^4-q^4$ and $2m$ must be the same as the ratio between a square and another square. Naturally, Diophantus does not tell us how he arrived at this statement, his system of algebraic notation being totally inadequate for this task. Also noteworthy seem to me both the lack of precision in the argument which Diophantus gives for the insolubility of $20x^2=3$ (``if only both coefficients would have been squares ...''), and the speed with which he moves from this argument to the claim that the problem is now reduced (\textgreek{ἀπάγεται}) to finding $p,q,m$ as described above: he evidently assumes the reader quite adept at algebra, not to mention familiar with the idea of ``backtracking'' after an initial substitution has failed.} He lets $p$ be indeterminate for the time being, and chooses the values $q=2$ and $m=p^2+4$. Now he has
$$
\frac{m^2-p^4-q^4}{2m}=\frac{(p^2+4)^2-p^4-2^4}{2m}=\frac{4p^2}{p^2+4},
$$
which means that $p^2+4$ needs to be made a square, and he puts $p^2+4=(p+1)^2$ which gives $p=\dfrac{3}{2}$. This gives $m=\dfrac{25}{4}$. On scaling $p$ and $q$ by $2$ and $m$ by $4$, he gets $p=3$, $q=4$, $m=25$. Then Diophantus goes back to his original equation, which now reads
$$
x^4+81+256 = (x^2 - 25)^2 ~~ \Longrightarrow ~~ x^4 + 337 =x^4-50x^2+625,
$$ 
which gives $50x^2 = 288$, so that $x = \dfrac{12}{5}$, which solves the problem.

\subsection{The geometry of the solution to Problem V.29}

Diophantus' problem asks for rational numbers $x,y,$ and $z$ such that there exists a rational number $w$, such that 
\begin{equation}
\label{V29a}
x^4+y^4+z^4=w^2.
\end{equation}

We shall consider the surface $X$ defined by \eqref{V29a} in the weighted projective space $\mathbb{P}(1,1,1,2)$, where $x,y,z$ have weight $1$ and $w$ has weight $2$. It is a del Pezzo surface $X$ of degree 2: since it is a double cover of $\mathbb{P}^2$ ramified in a smooth quartic, one may apply the Riemann--Hurwitz formula for generically finite morphisms (see \cite[p.~349]{reid}). 

On $X$, Diophantus constructs a rational curve as follows. He applies the substitution
$$
x^4-p^4-q^4=(x^2-m)^2
$$
and then asks which $p,q,m$ allow for a rational value of $x$. What this comes down to in geometric terms is that he  replaces the surface $X$ by the birationally equivalent surface
$$
X' \colon m^2-p^4-q^4 = 2mr^2,
$$
where again $X'$ is most conveniently described as a surface in weighted projective space $\mathbb{P}(1,1,1,2)$ with coordinates $m,p,q,r$, where $p,q,r$ have weight $1$ and $m$ has weight $2$. 

A birational map $f \colon X' \dashrightarrow X$ is given by
\begin{equation}
\label{V29bir}
(p,q,r,m) \mapsto (r,p,q,r^2+m),
\end{equation}
and it is easily verified that it is indeed birational. Following Diophantus' cue, we can construct a rational curve $C$ on $X'$ defined by the equation
$$
m = p^2+q^2.
$$
This allows us to eliminate $m$, which means projecting birationally to the projective plane $\mathbb{P}^2$ with coordinates $p,q,r$, to obtain the curve 
$$
C' \colon 2p^2q^2 = 2(p^2+q^2)r^2.
$$
But on $C'$, we can write 
$$
p^2+q^2=\left(\dfrac{pq}{r}\right)^2,
$$
showing that $C'$ is birationally equivalent to the unit circle. Hence the same holds for $C$. The image of $C$ is a curve on $X$ that is birationally equivalent to $\mathbb{P}^1$, since it is not contracted by $f$, as is easily seen from \eqref{V29bir}.

\section{Appendix: singular models of K3 surfaces}
\label{singK3}

In this short appendix, we will give a proof of Theorem \ref{singularK3},  which was used in Sections \ref{secIII17} and \ref{secIV18}. It is  by no means a new result, but since we have been unable to locate a proof in the literature, we have decided to include one here. I would like to thank Francesco Polizzi for references to the literature and some pointers regarding the proof.

We let $k$ be a field of characteristic $0$. With a \defi{surface over $k$}, we mean a scheme of dimension $2$ over $\operatorname{Spec} k$. We recall that a \defi{K3 surface} over $k$ is a regular, projective, and geometrically connected surface $Y$ over $k$ such that the canonical bundle $\omega_Y$ of $Y$ is trivial and the sheaf cohomology group $\operatorname{H}^1(Y,\mathscr{O}_Y)$ vanishes.

The main result of this appendix, which is well-known to experts, describes certain (possibly mildly singular) complete intersections which are birationally equivalent to K3 surfaces.

\begin{theorem}\label{singularK3}
Assume that $X$ is a surface over $k$ of one of the following three types:
\begin{itemize}
\item[(i)] a quartic surface in $\mathbb{P}^3_k$,
\item[(ii)] an intersection of a cubic and a quadric hypersurface in $\mathbb{P}^4_k$,
\item[(iii)] an intersection of three quadrics in $\mathbb{P}^5_k$.
\end{itemize} 
Furthermore, assume that all singularities of $X$ are rational double points. Then the minimal regular model of $X$ is a K3 surface.
\end{theorem}

Note that we do not explicitly require $X$ to be geometrically connected or geometrically integral; as we will see, this will follow from the other assumptions. For the definition of \defi{rational double points}, also called simple surface singularities, or canonical surface singularities, or Du Val singularities, see \cite{reid}. Intuitively speaking, they are the ``mildest'' kind of surface singularities. In particular, they are normal (e.g., because they are analytically isomorphic to hypersurface singularities of codimension $\geq 2$, see \cite[Corollary 8.2.24]{liu}).

\begin{proof}
We will show the existence of a birational morphism $$f \colon Y \rightarrow X$$ such that $Y$ is a K3 surface. It then follows that $Y$ is the unique surface with this property, by uniqueness of minimal regular models of surfaces of Kodaira dimension $\geq 0$ (see e.g. \cite{beauville}).

Let $i \colon X \rightarrow \mathbb{P}^n$ be a projective embedding of $X$ as a complete intersection of $n-2$ hypersurfaces, of degrees $d_1,\ldots,d_{n-2}$, with
\begin{itemize}
\item[(i)] $n=3$ and $d_1=4$, or
\item[(ii)] $n=4$ and $d_1=3, d_2=2$, or
\item[(iii)] $n=5$ and $d_1=2, d_2=2,d_3=2$.
\end{itemize} 

First, we show that the canonical sheaf $\omega_X$ of $X$ is trivial. We recall (cf.~\cite[Definition 6.4.7]{liu}) that one way of defining $\omega_X$ is
\begin{equation}
\label{omegaYdef}
\omega_X = \det i^\ast (\mathcal{I}/\mathcal{I}^2)^\vee  \otimes i^\ast \omega_{\mathbb{P}^n},
\end{equation}
where $\mathcal{I} \subset \mathscr{O}_{\mathbb{P}^n}$ is the ideal sheaf of $X$, and $\omega_{\mathbb{P}^n}$ denotes the canonical sheaf on $\mathbb{P}^n_k$. (See \cite[p.~349]{reid} for two equivalent definitions. We note that the more familiar definition of $\omega_X$ as the determinant of the sheaf of K\"ahler differentials is in general not valid when $X$ is non-regular, so we cannot use it here.) Proceeding as in the proof of \cite[Lemma IV.11]{beauville}, we consider the following surjective map of $\mathscr{O}_{\mathbb{P}^n}$-modules
$$
\mathscr{O}_{\mathbb{P}^n}(-d_1) \oplus \cdots \oplus \mathscr{O}_{\mathbb{P}^n}(-d_{n-2}) \twoheadrightarrow \mathcal{I}/\mathcal{I}^2.
$$
Applying the right-exact functor $i^\ast$ we get the surjection
$$
\mathscr{O}_{X}(-d_1) \oplus \cdots \oplus \mathscr{O}_{X}(-d_{n-2}) \twoheadrightarrow i^\ast ( \mathcal{I}/\mathcal{I}^2 )
$$
which is an isomorphism since both sides are locally free of rank $r-2$ (to see that this holds for $i^\ast(\mathcal{I}/\mathcal{I}^2)$, one uses that $X$ is a complete intersection; see \cite[Example 6.3.5, Corollary 6.3.8]{liu}). Taking determinants, we obtain $\det i^\ast(\mathcal{I}/\mathcal{I}^2 ) = \mathscr{O}_X (-d_1-\ldots -d_{n-2}) $. Substituting this into \eqref{omegaYdef}, and using that $\omega_{\mathbb{P}^n}$ is isomorphic to $\mathscr{O}_{\mathbb{P}^n}(-n-1)$, we get
\begin{equation}
\label{omegaYexpression}
\omega_X = \mathscr{O}_X (d_1+\ldots +d_{n-2}) \otimes i^\ast \omega_{\mathbb{P}^n} = \mathscr{O}_X.
\end{equation}
This proves the first claim.  

We now define $f \colon Y \to X$ as the minimal resolution of singularities of $X$. Since the singularities of $X$ are rational double points, we have $f^\ast \omega_X = \omega_Y$ (see \cite[p.~347]{reid}). We conclude that 
$$
\omega_Y = f^\ast \omega_X =  f^\ast \mathscr{O}_X = \mathscr{O}_Y. 
$$  

Our next claim is that $\operatorname{H}^1(X,\mathscr{O}_X)$ vanishes. (Our proof of this follows that of \cite[Lemma VIII.9]{beauville}. We reproduce the proof since we will use it for something else as well.) One reduces to the following, more general claim: {\it if $X_r$ is a complete intersection of $r$ hypersurfaces in $\mathbb{P}^n_k$, then 
$$
\operatorname{H}^j(X_r,\mathscr{O}_{X_r}(m))=0
$$
for all integers $0 < j < n-r$ and all integers $m$.} Applying induction, with $X_0 = \mathbb{P}^n_k$ as the (trivial) base case, we may assume that the result holds for any complete intersection $X_{r-1}$ of $r-1$ hypersurfaces in $\mathbb{P}^n$, and that $X_r \to X_{r-1}$ is a closed immersion defined by a homogeneous form of degree $d$. The claim now immediately follows from considering the long exact cohomology sequence associated to
\begin{equation}
\label{cohomology}
0 \rightarrow \mathscr{O}_{X_{r-1}}(m-d) \rightarrow \mathscr{O}_{X_{r-1}}(m) \rightarrow \mathscr{O}_{X_{r}}(m)  \rightarrow 0.
\end{equation}

It follows from the same long exact sequence that $X$ is geometrically connected. Indeed, taking $m=0$ in \eqref{cohomology}, the first part of the sequence is as follows
\begin{align*}
0 \rightarrow \operatorname{H}^0(X_{r-1},\mathscr{O}_{X_{r-1}}(-d)) & \rightarrow \operatorname{H}^0(X_{r-1},\mathscr{O}_{X_{r-1}}) \\ \rightarrow & \operatorname{H}^0(X_{r},\mathscr{O}_{X_r}) \rightarrow \operatorname{H}^1(X_{r-1},\mathscr{O}_{X_{r-1}}(-d)),
\end{align*}
from which we deduce $\operatorname{H}^0(X_r,\mathscr{O}_{X_r})=\operatorname{H}^0(X_{r-1},\mathscr{O}_{X_{r-1}})$ for all $r \geq 1$, so that $\operatorname{H}^0(X_r,\mathscr{O}_{X_r})=\operatorname{H}^0(\mathbb{P}^n_k,\mathscr{O}_{\mathbb{P}^n})=k$, proving that $X_r$ is geometrically connected. Hence $X$ is geometrically connected, and therefore also geometrically integral: indeed, if $\overline{x}$ is a point on $\overline{X} := X \times_k \overline{k}$, the local ring $\mathscr{O}_{\overline{X},\overline{x}}$ at $\overline{x}$ is either regular or the local ring of a rational double point; in either case, it is clear that $\mathscr{O}_{\overline{X},\overline{x}}$ is an integral domain, and this implies that $\overline{x}$ is contained in only one irreducible component of $\overline{X}$. Since $Y$ arises from $X$ by repeated blowing-up of singularities, it is also geometrically integral.

We now claim that $\operatorname{H}^1(Y,\mathscr{O}_Y)$ also vanishes. Note that, by the fact that $f \colon Y \rightarrow X$ is a projective birational morphism between geometrically integral varieties, and $X$ is normal, we have $f_\ast \mathscr{O}_Y = \mathscr{O}_X$ (compare the proof of \cite[Corollary III.11.4]{hartshorne}). We use the Leray spectral sequence associated to the map $f$ and the sheaf $\mathscr{O}_Y$:
$$
E^2_{p,q} = \operatorname{H}^p(X,\operatorname{R}^q f_\ast \mathscr{O}_Y) \Rightarrow \operatorname{H}^{p+q}(Y,\mathscr{O}_Y).
$$
Since the singularities of $X$ are rational double points, the higher direct image $\operatorname{R}^1 f_\ast \mathscr{O}_Y$ vanishes \cite[Section 3]{reid}. Hence $E^2_{0,1}=0$, which implies $E^\infty_{0,1}=0$, and $E^\infty_{1,0}=E^2_{1,0}=\operatorname{H}^1(X, f_\ast \mathscr{O}_Y)$. We deduce
$$
\operatorname{H}^1(Y,\mathscr{O}_Y) = \operatorname{H}^1(X,f_\ast \mathscr{O}_Y) = \operatorname{H}^1(X,\mathscr{O}_X)=0.
$$

Finally, $Y$ is projective and regular by construction. Furthermore, we saw that $Y$ is geometrically connected, and it is of dimension $2$ since it is birational to $X$. Hence $Y$ is a K3 surface.
\end{proof}

\bibliography{diophantus}
\bibliographystyle{plain}
\end{document}